\documentclass{amsart}

\usepackage{amssymb, amsmath,mathtools}
\usepackage{mathrsfs,mathabx}
\usepackage{amscd}
\usepackage{verbatim}
\usepackage[dvipsnames]{xcolor}
\usepackage{subfig}

\usepackage{enumerate}
\usepackage{stmaryrd}

\usepackage[colorlinks,linkcolor={blue},citecolor={blue},urlcolor={purple},]{hyperref}

\usepackage[colorinlistoftodos,prependcaption,textsize=tiny]{todonotes}

\theoremstyle{plain}
\newtheorem{theorem}{Theorem}[section]
\theoremstyle{remark}
\newtheorem{remark}[theorem]{Remark}
\newtheorem{example}[theorem]{Example}
\theoremstyle{plain}

\newtheorem{lemma}[theorem]{Lemma}
\newtheorem{proposition}[theorem]{Proposition}

\newtheorem{definition}[theorem]{Definition}

\newtheorem{assumption}[theorem]{Assumption}

\numberwithin{equation}{section}

\usepackage{a4wide}
\def\N{{\mathbb N}}
\def\Z{{\mathbb Z}}

\def\R{{\mathbb R}}

\newcommand{\E}{{\mathbf E}}
\renewcommand{\P}{{\mathbf P}}

\newcommand{\g}{\gamma}

\newcommand{\om}{\omega}

\newcommand{\dd}{\mathrm{d}}
\newcommand{\T}{\mathbb{T}}

\newcommand{\Lip}{\mathrm{Lip}}
\newcommand{\BiLip}{\mathrm{BiLip}}
\newcommand{\Diff}{\mathrm{Diff}}
\newcommand{\Prob}{\mathcal{P}}
\newcommand{\X}{\mathfrak{X}}
\newcommand{\dif}{\mathrm{d}}

\newcommand{\Wasl}{\mathcal{W}^{\Lip}}

\newcommand{\eps}{\varepsilon}

\DeclareMathOperator{\Law}{Law}

\DeclareMathOperator{\tr}{tr}

\newcommand{\into}{\hookrightarrow}
\DeclarePairedDelimiter\ceil{\lceil}{\rceil}
\DeclarePairedDelimiter\floor{\lfloor}{\rfloor}
\newcommand{\dist}{\mathrm{d}}
\newcommand{\D}{\mathscr{D}}

\begin{document}

\author{Antonio Agresti}
\address{Department of Mathematics Guido Castelnuovo\\ Sapienza University of Rome\\
  P.le Aldo Moro 5\\ 00185 Rome\\ Italy}
\email{antonio.agresti@uniroma1.it}

\author{William Cooperman}
\address{Department of Mathematics\\ New York University \\ 251 Mercer Street\\ New York\\ NY 10012\\ USA}
\email{bill@nyu.edu}

\author{Patrick Flynn}
\address{Department of Applied \& Computational Mathematics\\ Yale University\\ 51 Prospect Street \\ New Haven \\ CT 06511\\  USA }
\email{patrick.t.flynn@yale.edu}

\date\today

\title[Stability of negative-regularity mixing and applications]{Stability of negative-regularity mixing and\\ applications to stochastic fluids}




\begin{abstract}
We consider mixing for a scalar advected by random vector fields, such as solutions to stochastic fluid equations.
Our main abstract result is that contraction in expectation of a negative Sobolev norm of the scalar (negative-regularity mixing) is stable under perturbations of the joint scalar--velocity process.
We apply the result in three settings: mixing for a hyperviscous Boussinesq-type system where the scalar actively drives the velocity, mixing for a passive scalar advected by two-dimensional Navier--Stokes driven by smooth-in-time stochastic forcing, and a black-box proof that negative-regularity mixing implies enhanced dissipation.
\end{abstract}
\maketitle

\section{Introduction}


This paper concerns the long-time behavior of a scalar $\theta_t$ transported by a
divergence-free random vector field $u_t$, both defined on a closed manifold $M$:
\begin{equation}
  \label{eq:advection}
    \partial_t \theta_t +  u_t\cdot  \nabla \theta_t = 0 \quad  \text{  on }\   M\times [0,\infty),
\end{equation}
We seek to identify sufficient conditions for \textit{scalar mixing}: that is, whether or not $\theta_t$ converges weakly to a constant as $t \to \infty$. In the present paper, we restrict our study to \textit{exponential} mixing.
To be precise, assuming $\int \theta_0\dif  x = 0$, one may hope to show a \textit{quenched mixing estimate} of the form
\begin{equation}
  \label{eq:mixing-estimate}
  \|\theta_t \|_{H^{-1}(M)} \leq C\|\theta_0\|_{H^{1}(M)}e^{-\gamma t},
\end{equation}
where $\gamma > 0$, $C = C(u) > 0$ is a random constant, and $H^s(M)$ denotes the Sobolev space with regularity $s$. When $u$ is autonomous and deterministic, Dolgopyat proved~\eqref{eq:mixing-estimate} for a large class of uniformly hyperbolic flows.
Later,~\eqref{eq:mixing-estimate} was shown to hold for a large class of random processes $(u_t)_{t\geq 0}$, including many fundamental models in stochastic fluid mechanics~\cite{BBPS22}. We refer to Section~\ref{sec:background} below for a brief history of fluid mixing in general.

The recent work~\cite{BFPS24} (see also~\cite{dewitt2025conservative}) proves  \textit{negative regularity annealed mixing}\footnote{
The term ``annealed'' refers to decay in expectation, while ``quenched'' refers to almost sure decay.} for a large class of random volume preserving flows:
\begin{equation}
  \label{eq:negreg-mixing-estimate}
  \E \|\theta_t\|_{H^{-\delta}(M)} \leq C\|\theta_0\|_{H^{-\delta}(M)}e^{-\gamma t},
\end{equation}
where $\delta,\gamma> 0$ are small, and $C$ is a deterministic constant that depends only on the domain $M$ and the random process $(u_t)_{t\geq 0}$.

In contrast to~\eqref{eq:mixing-estimate}, the negative-regularity estimate~\eqref{eq:negreg-mixing-estimate} shows decay of $\theta_t$ in the negative regularity Sobolev space $H^{-\delta}(M)$ assuming only that $\theta_0 \in H^{-\delta}(M)$; mean-zero scalars decay with \textit{no loss of derivatives}. Since time-reversibility of~\eqref{eq:advection} precludes almost-sure contraction, uniform in $\theta_0$, in any fixed space, estimate~\eqref{eq:negreg-mixing-estimate} records decay of $\|\theta_t\|_{H^{-\delta}}$ in expectation only.
Despite the expectation,~\eqref{eq:negreg-mixing-estimate} is strictly stronger than~\eqref{eq:mixing-estimate}; the quenched $H^1(M) \to H^{-1}(M)$ estimate can be recovered (at the expense of taking $\gamma$ smaller and $C$ larger) through a Borel--Cantelli argument.

\subsection{Illustration of the main results}
In this paper, we show that such negative regularity mixing estimates~\eqref{eq:negreg-mixing-estimate} are stable under various perturbations of the joint process $(\theta_t,u_t)$ under general assumptions (see Assumptions~\ref{ass:negative_regularity_mixing} and~\ref{ass:super_lyapunov}). 
We separate our stability result into two main abstract theorems: 
\begin{itemize}
    \item Theorem~\ref{thm:mainStability} shows that for some fixed finite time (for instance, $t = 1$), the estimate~\eqref{eq:negreg-mixing-estimate} is stable under perturbations of the process $((\theta_t,u_t))_{t\geq 0}$.
    \item Theorem~\ref{thm:shortTimeToLong} proves that a finite time estimate may be iterated to prove that~\eqref{eq:negreg-mixing-estimate} holds for arbitrarily large times as well.
\end{itemize}
We apply the abstract results in three separate settings, which we introduce now.

\subsubsection{Enhanced dissipation} In Section~\ref{sec:enhanced-dissipation} we perturb~\eqref{eq:advection} with a small diffusive term:
\begin{equation}
  \label{eq:advection-diffusion}
  \begin{cases}
    \partial_t \theta_t^\kappa + u_t \cdot \nabla \theta^\kappa_t = \kappa \Delta \theta^\kappa_t &\quad \text{$t > 0$,}\\
    \theta^\kappa_0 = \theta_0 &\quad \text{$t = 0$,}
  \end{cases}
\end{equation}
where $\kappa > 0$. In Theorem~\ref{thm:enhanced-dissipation},
we show that the negative regularity mixing estimate~\eqref{eq:negreg-mixing-estimate} is stable under perturbation of the velocity by small Brownian translations, and by taking an expectation we conclude that~\eqref{eq:negreg-mixing-estimate} holds for~\eqref{eq:advection-diffusion} with $\kappa$-independent constants. As a consequence, we prove that flows exhibiting negative regularity mixing enhance dissipation:
\begin{align}
\|\theta^\kappa_t\|_{L^2} \leq D\kappa^{-{\delta/2}} e^{-\gamma t}e^{\chi V(u_0)}\|\theta_0\|_{L^2}, \qquad t\geq 0,
\end{align}
where $\delta$ is as in~\eqref{eq:negreg-mixing-estimate}. This result complements similar works on enhanced dissipation of exponentially mixing vector fields~\cite{ConstantinEtAl2008,CZDT20}, see Subsection~\ref{sss:literature} for further references.

\subsubsection{Mixing estimates for thermally driven weakly active scalars} In Section~\ref{s:thermal_flows}, as a second application, we show that negative regularity mixing occurs for a fluid coupled to a weakly active scalar. We call the scalar \textit{weakly active} because it only slightly influences the velocity; our results only apply with $\varepsilon$ in~\eqref{eq:eps-NSE-intro} taken small.
We prove~\eqref{eq:negreg-mixing-estimate} for incompressible Navier--Stokes on $\T^2 = \R^2/\Z^2$ with stochastic forcing and a small coupling term involving the scalar:
\begin{equation}\label{eq:eps-NSE-intro}
  \begin{cases}
    u^\varepsilon = \nabla^\perp \Delta^{-1} \omega^\varepsilon\\
    \partial_t \omega^\varepsilon + u^\varepsilon \cdot \nabla \omega^\varepsilon +  \Delta^2 \omega^\varepsilon = \dot{F} - \varepsilon \partial_{1}\theta^\varepsilon &\quad \text{for } t > 0,\\
    \partial_t \theta^\varepsilon + u^\varepsilon \cdot \nabla \theta^\varepsilon = 0 &\quad \text{for } t > 0,\\
    \omega^\varepsilon = \omega_0, \quad \theta^\varepsilon = \theta_0 &\quad \text{for } t = 0,
  \end{cases}
\end{equation}
where $\varepsilon > 0$ and the noise $\dot{F}$ is white in time and colored in space; see Section~\ref{s:thermal_flows} for the precise setup.

Since the scalar feeds back into the velocity, one must work with the evolution of the full scalar $\theta$ rather than tracking individual particles. Previous methods used to prove mixing (e.g.~\cite{DolgopyatKaloshinKoralov2004, BBPS22, BlumenthalCotiZelatiGvalani2023}) rely on proving unique geometric ergodicity of the two-point process, which is not even well-defined in this case.

However, we show that the conclusions reached using such Lagrangian methods are stable: for sufficiently small $\varepsilon > 0$ and $\|\theta_0\|_{L^\infty} \leq 1$, we have mixing
\begin{equation}
  \E \|\theta^\varepsilon_t\|_{H^{-\delta}(M)} \leq C\|\theta_0\|_{H^{-\delta}(M)}e^{-\gamma t}
\end{equation}
with constants $C, \gamma > 0$ independent of $\varepsilon$; see Theorem~\ref{t:hNSE-mixing}.


\subsubsection{Mixing for 2D Navier-Stokes with smooth random forcing} As our third application, in Section~\ref{sec:regularized}, we consider the 2D Navier-Stokes equation with a random smooth body force:
\begin{equation}
  \label{eq:2dnse-smooth-noise_intro}
  \begin{cases}
    \partial_t \omega^\varepsilon + u^\varepsilon \cdot \nabla \omega^\varepsilon - \nu \Delta \omega^\varepsilon = \xi^\varepsilon &\qquad \text{ for }x \in \mathbb T^2, t>0,\\
    \omega^\varepsilon= \omega_0, &\qquad\text{ for } x \in \mathbb T^2, t = 0.
  \end{cases}
\end{equation}
Above, $\xi^\eps$ is taken to be a smooth mollification in time of colored-in-space, white-in-time noise, where $\varepsilon > 0$ denotes the mollification radius; see~\eqref{eq:F-eps-def}.
To the best of our knowledge, it remains unknown whether this model, as well as other similar infinite-dimensional fluid models with smooth random forcing, exhibit the same unique ergodicity properties as their counterparts with white-in-time forcing. Nevertheless, we show in Theorem~\ref{thm:reg-noise-mixing} that the vector field $u^\eps =\nabla^\perp \Delta^{-1}\omega^\eps$ exhibits negative regularity scalar mixing as in~\eqref{eq:negreg-mixing-estimate} when $\varepsilon > 0$ is sufficiently small.
As mentioned above, the latter also implies the quenched estimate~\eqref{eq:mixing-estimate}.

\subsection{Background}\label{sec:background}\enlargethispage{4pt}

The theory of mixing is a broad subject at the intersection of statistical physics, fluid mechanics, dynamical systems, and probability. While systems such as~\eqref{eq:advection} are formally time-reversible, they often exhibit irreversible behavior after statistical averaging or examining the system over very long times. In this sense, mixing is one manifestation of Boltzmann's $H$-theorem, or the second law of thermodynamics. Indeed, Gibbs~\cite{Gibbs1902} used the example of ink in water\footnote{more precisely, a ``liquid to contain a certain amount of coloring matter''} as an analogy to explain the emergence of entropy in Hamiltonian systems.

\subsubsection{Dynamical systems} Contemporaneously with the development of ergodic theory~\cite{BirkhoffKoopman1932,Birkhoff1931,Koopman1931,KoopmanVonNeumann1932,vonNeumann1932}, Hopf gave a mathematical definition of Gibbs' notion of mixing as $L^2$ weak convergence (or decorrelation)~\cite{Hopf1932Transitivity,Hopf1932Gibbs}.
Since then, the theory of mixing has been developed from several perspectives. Perhaps the most fruitful of these approaches is that of hyperbolic dynamics, beginning with the work of Hopf and Hedlund~\cite{Hedlund1936,Hopf1939} on geodesic flows of hyperbolic surfaces. A mathematical theory of entropy was later developed in the works of~\cite{Kolmogorov1958,Kolmogorov1959,Sinai1959}, after which
hyperbolic dynamics proved to be fertile ground for such theories throughout the 20th century~\cite{Anosov1969,Bowen1975,HoweMoore1979,Pollicott1984,Pollicott1985,Ruelle1968,Ratner1987,Ruelle1983,Sinai1972}.
For uniformly hyperbolic \textit{flows} specifically, this theory culminated in the aforementioned work of Dolgopyat~\cite{Dolgopyat1998}, which proved quantitative exponential mixing for a general class of such systems. Numerous works have extended, generalized or otherwise further explored exponential mixing in other settings, e.g.~\cite{Pollicott1999,Stoyanov2011,Tsujii2018}. In particular, we highlight the works~\cite{BaladiTsujii2007,BlankKellerLiverani2002,Liverani2004,ButterleyWar2020,Fried1995,Rugh1992,Rugh1996} on anisotropic Banach spaces, as well as the semiclassical treatments of such spaces~\cite{FaureRoy2006,FaureRoySjostrand2008,FaureSjostrand2011,FaureTsujii2013,FaureTsujii2017,FaureTsujii2023,FaureTsujii2024,DyatlovZworski2015,DyatlovFaureGuillarmou2015,DyatlovGuillarmou2016,DyatlovZworski2016}. This approach was useful in the works on negative regularity mixing preceding this paper~\cite{BFPS24,DeWittDolgopyat2024,dewitt2025conservative,dewitt2026ergodicity}.

\subsubsection{Deterministic and stochastic fluids}
In parallel with the aforementioned mathematical works on dynamical systems throughout the 20th century, the problem of passive scalar mixing was actively studied by the fluid dynamics community~\cite{taylor1922,Obukhov1949,Eckart1948,Corrsin1951,Kraichnan1968,Pierrehumbert1994,Batchelor1959,Aref1984,Ottino1989}. A major motivation for studying  passive scalar mixing, and indeed one of the central problems of fluid dynamics, is the development of a theory of turbulence. In particular, such a theory would explain the phenomenological theory put forth by Kolmogorov in his famous 1941 papers~\cite{kolmogorov1941local,kolmogorov1941dissipation}. While such a theory for nonlinear fluid equations remains elusive, the analysis of their passive scalars is relatively more tractable, and yields some indirect understanding of turbulent transport.

 Around the turn of the millennium, major progress was made on the mathematical theory of stochastic fluid equations, such as the existence and uniqueness of ergodic measures for 2D Navier-Stokes equations with various stochastic forcings~\cite{flandoli1994dissipativity,flandoli1995martingale,FlandoliMaslowski1995,ferrario1997ergodic,BricmontKupiainenLefevere2002,KuksinShirikyan2002,KuksinShirikyan2004,HairerMattingly2006,HairerMattingly2008}. Around the same period, the mathematical fluids community became interested in scalar mixing for deterministic flows~\cite{Bressan2003,CrippaDeLellis2008,ConstantinEtAl2008}.

 Following these works, the mathematical theory of mixing for deterministic and stochastic fluids has been an area of intense research activity. For instance, in the deterministic setting, there have been numerous works on the (sub-)optimality of mixing estimates, constructions   of exponentially mixing vector fields, or the relationship between exponential mixing and enhanced dissipation~\cite{Seis2013,IyerKiselevXu2014,ElgindiZlatos2019,FengIyer2019,CotiZelatiDrivas2021,Seis2023Enhanced,ElgindiLissMattingly2025,cooperman2025exponentially}. See also~\cite{zelati2024mixing} for a recent survey of  the subject. We also note that similar progress has been made on the  closely related problems of anomalous dissipation ~\cite{DrivasElgindiIyerJeong2022,ColomboCrippaSorella2023,BurczakSzekelyhidiWu2023,ElgindiLiss2024,ArmstrongVicol2025,HessChildsRowan2025},
and nonlinear phase mixing~\cite{MouhotVillani2011,BedrossianMasmoudi2015,BedrossianMasmoudiVicol2016,BedrossianCotiZelati2017}.

In the stochastic fluids setting, almost-sure exponential mixing was established for various stochastic fluid equations (in particular, 2D stochastic Navier-Stokes) in the work~\cite{BBPS22}, see also~\cite{BedrossianBlumenthalPunshonSmith2022Lagrangian,BedrossianBlumenthalPunshonSmith2021}. This result was then extended to the case of stochastic Navier-Stokes with degenerate forcing in~\cite{CR26}. Similar results on mixing, enhanced dissipation or anomalous dissipation have also been proved for a variety of random fluid flows~\cite{BlumenthalCotiZelatiGvalani2023,Cooperman2023,CoopermanIyerSon2025,NavarroFernandezSeis2025,CotiZelatiNavarroFernandez2026}, and  fluids with transport noise~\cite{flandoli2021mixing,GessYaroslavtsev,Agresti2025Transport,Rowan2024Kraichnan}. Finally, we mention recent progress on spectral properties of turbulent advection diffusion equations, such as  the Batchelor spectrum~\cite{BedrossianBlumenthalPunshonSmith2022Batchelor,BlumenthalHuynh2026,cooperman2026fourier},  intermittency~\cite{ColomboCrippaSorella2023,ArmstrongVicol2025,HessChildsRowan2025Intermittency}, and the Obukhov--Corrsin spectrum~\cite{rowan2025obukhov}.

\subsubsection{Relationship with other works and future directions}
\label{sss:literature}
As mentioned previously, the main abstract results of this paper, Theorems~\ref{thm:mainStability} and~\ref{thm:shortTimeToLong},  show that the negative regularity mixing results of~\cite{BFPS24} are stable under fairly general perturbations of the underlying stochastic process. Like~\cite{BFPS24}, we assume the existence of a certain well-behaved super-Lyapunov function $V(u_t)$ in our assumptions, which appears as a weight in our estimates. Propagating such weights is necessary in situations where the state space of $u_t$ is noncompact, and the desired mixing bounds are nonuniform in $u_t$. In particular, such super-Lyapunov functions are necessary to prove mixing by vector fields defined by solutions to dissipative stochastic PDEs, such as 2D stochastic Navier-Stokes, as well as other fluid equations considered in~\cite{BBPS22,BFPS24}.

While our abstract stability result allows us to prove negative regularity mixing for a variety of systems, the stability of mixing under perturbations is not without precedent. Indeed, our result could be viewed as an analogue of the spectral stability result of~\cite{keller1999stability}, although the abstract setting and the nature of the result are different. See also~\cite{TsujiiZhang2026} for a very recent result in the Anosov setting.

We now compare the various applications to other works. First, Theorem~\ref{thm:enhanced-dissipation} is in the same spirit as the main results of~\cite{ConstantinEtAl2008,CZDT20}, which prove that mixing of the undamped equation implies enhanced dissipation of the damped equation. In contrast to~\cite{ConstantinEtAl2008}, we make no use of the RAGE Theorem. Moreover, since we assume a negative regularity mixing estimate~\eqref{eq:negreg-mixing-estimate}, and consider random vector fields $u^\eps$, our theorem applies in a different setting.

In Section~\ref{s:thermal_flows}, we establish negative regularity mixing for active scalar systems. This appears to be entirely new to the stochastic PDE literature, as the only other setting where mixing is known is for nonlinear transport equations near shear flows, i.e. Landau damping~\cite{MouhotVillani2011} and inviscid damping for 2D Euler~\cite{BedrossianMasmoudi2015}. It would be interesting to investigate a similar result for an inviscid stochastic fluid equation, although this goes well beyond the capabilities of our current perturbation theory. Indeed, even replacing hyperviscosity with ordinary viscosity in Section~\ref{s:thermal_flows} poses an interesting technical challenge.

Finally, in Section~\ref{sec:regularized}, we prove negative regularity mixing for 2D stochastic Navier-Stokes with regularized stochastic forcing. While unique ergodicity and mixing for Navier-Stokes for a smooth bounded random forcing has been established~\cite{Shirikyan2015, KuksinNersesyanShirikyan2020}, obtaining the same results where the forcing is a mollification-in-time of white noise remains open. Without proving anything about the velocity process as an intermediate step, we are able to conclude that passive scalars mix exponentially in the smooth-in-time setting.




\medskip

\noindent{\bf Notation.}
We write $C, c > 0$ throughout to refer to large and small constants (respectively) that may change from line to line. We denote the unit torus by $\mathbb T = \mathbb R/\mathbb Z \simeq S^1$. We use $\mathbf E$ and $\mathbf P$ to denote expectation and probability measure respectively with respect to the ambient probability space when applicable.

\medskip

\noindent{\bf Acknowledgments.} The authors would like to thank Jacob Bedrossian, Alex Blumenthal, Sam Punshon-Smith and Keefer Rowan for useful discussions and suggestions. The original concept for this work was developed by the authors during the workshop ``Random Dynamical Systems, PDEs, and Stochastic Analysis" at the Brin Mathematics Research Center at the University of Maryland.
We thank the organizers and staff for their hospitality and support.
Antonio Agresti is a member of GNAMPA (INdAM) and acknowledges support from INdAM through the GNAMPA 2026 project ``Fluidodinamica stocastica: irregolarità, trasporto e fenomeni di regolarizzazione''.
William Cooperman was partially supported by National Science Foundation grant DMS-2303355.
Patrick Flynn gratefully
acknowledges support from the Simons Laufer Mathematical Sciences Institute,
where part of the research was conducted during the program ``Kinetic Theory: Novel Statistical, Stochastic
and Analytical Methods."

\medskip

\noindent{\bf AI statement.} The core ideas of this paper were developed by the authors in the summer of 2025. In recent months, we used publicly available models from Anthropic and OpenAI to assist in background research, sharpening results, proving  remaining technical lemmas, checking for errors, and proofreading in general. All content is written by the authors, and all AI suggested revisions have been carefully reviewed by the authors before implementing. We take full responsibility for any mistakes.
\section{Preliminaries}

\subsection{Calculus on manifolds}
Throughout this manuscript, we assume that $M$ is a closed Riemannian manifold of dimension $d\geq 2$. We denote the set of smooth diffeomorphisms from $M$ to $M$ by $\mathrm{Diff}(M)$, and denote the set of smooth vector fields on $M$ by $\X(M)$. Below we define the set of finite regularity diffeomorphisms, and their associated metrics.

\begin{definition}
  Fix a smooth isometric embedding $\iota \colon M \into \R^L$ for some large $L \in \N$.

  For $k \geq 0$, we will use the $C^k$-distance on diffeomorphisms given by
  \begin{equation}
    \label{eq:diffDistance}
    d^k(\varphi, \psi) \coloneq d^k_\to(\varphi, \psi) + d^k_\to(\varphi^{-1}, \psi^{-1}),
  \end{equation}
  where
  \begin{equation}
    \label{eq:diffDistanceOneWay}
    d^k_\to(\varphi, \psi) \coloneq \|\iota \circ \varphi - \iota \circ \psi\|_{C^k}.
  \end{equation}
  Write $\Diff^k(M)$ to denote the closure of $\Diff(M)$ under the metric $d^k$. Similarly, we denote $\X^k(M)$ to be the closure of $\X(M)$ with respect to  the norm $u \mapsto \|d\iota\,u\|_{C^k}$.
  Write $I\in \Diff(M)$ to denote the identity map $I \colon x \mapsto x$.
  \end{definition}

It will also be convenient to introduce Lipschitz norms. Given $\varphi \in \Diff^1(M)$, we write
$$
\BiLip(\varphi) \coloneq \Lip(\varphi) + \Lip(\varphi^{-1}),
$$
where $\Lip$ denotes the standard Lipschitz seminorm
$\Lip(\varphi) \coloneq \sup_{x\neq y} \frac{\dist(\varphi(x), \varphi(y))}{d(x,y)}$. One can check that
$$
\operatorname{Lip}(\varphi)
=
\sup_{x\in M}\|D\varphi_x\|_{T_x M\to T_{\varphi(x)} M}.
$$

Next, we record a key lemma relating the stability of flow maps on $M$ to their respective vector fields.

\begin{lemma}\label{lem:velocityToFlowMap}
  Fix $\alpha > 0$ and $k \in \mathbb N_0$ with $k+\alpha\geq 1$, and set $r=k+\alpha$.\par
  Let $u,v\in L^1_{\mathrm{loc}}([0,\infty);\X^r(M))$, and let $\varphi_t,\psi_t$ be their flow maps from time $0$ to time $t$.
  Then
  \begin{equation}
    \label{eq:velocityFlowEstimate}
    d^k(\varphi_t, \psi_t) \leq C \left(\int_0^t\|u_s - v_s\|_{L^\infty} \dif s\right)^{\frac{\alpha}{k+\alpha}} \exp\left(C \int_0^t \left(\|u_s\|_{C^{r}} + \|v_s\|_{C^{r}}\right) \dif s\right),
  \end{equation}
and
  \begin{equation}
    \label{eq:flowMapNormGrowth}
    d^{k+\alpha}_{\to}(I, \varphi_t) \leq C \exp\left(C \int_0^t \|u_s\|_{C^{k+\alpha}} \dif s\right).
  \end{equation}
\end{lemma}
\begin{proof}
  The estimate~\eqref{eq:velocityFlowEstimate} $k = 0$ follows by Gr\"onwall.
  When $k \geq 1$, it suffices to prove~\eqref{eq:flowMapNormGrowth}
  and deduce~\eqref{eq:velocityFlowEstimate} by interpolation (first replacing $d^k_\to$ with $d^k$ by symmetry).

  Working in a local coordinate patch, we take $k$ derivatives of $\dot{\varphi}_t(x) = u_t(\varphi_t(x))$ to obtain
  \begin{equation*}
    \partial_t \|\varphi_t\|_{C^{k+\alpha}} \leq C\left(\|u_t\|_{C^{k+\alpha}} \|\varphi_t\|_{C^{k-1+\alpha}}^k + \|u_t\|_{C^{1+\alpha}}\|\varphi_t\|_{C^{k+\alpha}}\right).
  \end{equation*}
  By induction on $k$ and Gr\"onwall, we conclude~\eqref{eq:flowMapNormGrowth}.
\end{proof}

For $\varphi\in\Diff(M)$, let $J_\varphi=|\det D \varphi(x)|$ denote its Riemannian Jacobian, where $D\varphi(x): T_x M \to T_{\varphi(x)} M$ is the standard differential.
For a measurable function $\theta:M\to\mathbb R$, we define
$$
(\varphi_*\theta)(x)
\coloneq
J_{\varphi^{-1}}(x)\theta(\varphi^{-1}(x)).
$$
For $\theta\in\D'(M)$ and $f\in C^\infty(M)$, define the pushforward by
$$
\langle \varphi_*\theta,f\rangle
=
\langle \theta,f\circ\varphi\rangle.
$$
Given a subspace of distributions $\mathscr{X}\subseteq \D'(M)$ (such as \ $L^2(M)$, or Sobolev spaces $H^s(M)$ defined below), we denote by
$$
\mathscr{X}_0 = \{\eta \in \mathscr X \ : \ \langle \eta, 1\rangle = 0\}
$$
the subspace of mean-zero distributions in $\mathscr{X}$.

\subsection{Sobolev spaces on manifolds}

Let $\nabla$ be the (weak) Riemannian gradient. We denote by $H^k(M)$ the set of all $f\in L^2(M)$ such that $\nabla^k f\in L^2(M;TM^{\otimes k})$, and define
$$
\| f\|_{H^k(M)}
\coloneq \|f\|_{L^2(M)}+\|\nabla^k f\|_{L^2(M;TM^{\otimes k})}.
$$
Moreover, for $\delta\in (0,1)$, we let
\begin{equation}
\label{eq:definition_fractional_Sob_spaces}
H^{-1}(M)\coloneq(H^1(M))^* \quad \text{ and }\quad H^{-\delta}(M)= [L^2(M),H^{-1}(M)]_\delta.
\end{equation}
where $[X,Y]_\delta$ denotes the complex interpolation space between two Banach spaces $X$ and $Y$ with the interpolation parameter $\delta$ (see~\cite[Appendix C.2]{Analysis1} for precise definitions).\footnote{We can also define $H^s$ for fractional $s >0$ this way as well.} In particular, the theory gives us an exact interpolation inequality,
\[
\|f\|_{H^{-\delta}(M)} \leq \|f\|_{L^2(M)}^{1-\delta} \|f\|_{H^{-1}(M)}^\delta, \qquad f \in L^2(M) \cap H^{-1}(M).
\]
We note that this is equivalent to the fractional Sobolev spaces on manifolds defined by localization in~\cite[Chapter 4, Subsection 3]{TayPDE1}, although the inequality above only holds up to a constant.

When the domain is unambiguous, we write $L^2 = L^2(M)$, $H^s=H^s(M)$, and so on.


\section{Setup and abstract results}
\label{sec:abstract_results}
\subsection{Random diffeomorphisms driven by Markov chains}
Once again, we let $M$ be a closed Riemannian manifold of dimension $d$. Given a measurable space $(X,\mathcal X)$, we shall use $\Prob(X)$ to denote the set of $\mathcal X$-measurable probability measures on $X$. If $X$ is a topological space, we take $\mathcal X$ to be the Borel $\sigma$-algebra unless otherwise specified.

We fix a measurable space $(Z,\mathcal Z)$, define $P_1 : Z \times \mathcal Z \to [0,1]$  to be  Markov transition kernel on $Z$. More precisely: first, we assume for each $A \in \mathcal Z$, the function $z\mapsto P_1(z,A)$ is $\mathcal Z$-measurable; second, $P_1(z) \in \mathcal P(Z)$ for each $z \in Z$.
We define the transition kernel at later times $
P_t \colon Z \to \Prob(Z),  t \in \mathbb N_0$ by recursively applying the semigroup property $P_{t+s}(z) = \int P_s(w) \dif P_t(z)(w)$ for each $z \in Z$ and $s,t \in \N_0$. Here, we set $P_0(z) = \delta_z$.

Recall that we  may then define joint transition probabilities of the Markov chain, and ultimately the measure on paths in $Z^{\mathbb N_0}$, as follows:
  given $t_1, \dots, t_n \in \N$ with $t_1 < t_2 < \dots < t_n$, define $P_{t_1, \dots, t_n} \colon Z \to \Prob(Z^n)$ recursively by
  \begin{align*}
    &\int_{Z^n} f(z_1, \dots, z_n) \dif P_{t_1, \dots, t_n}(z)(z_1, \dots, z_n)\\
    &\qquad \coloneq \int_Z \int_{Z^{n-1}} f(z_1, \dots, z_n) \dif P_{t_n - t_{n-1}}(z_{n-1})(z_n) \dif P_{t_1, \dots, t_{n-1}}(z)(z_1, \dots, z_{n-1}).
  \end{align*}
  For brevity, we let
  $$
  P_{[1,n]} \coloneq P_{1,2,\dots,n}.
  $$
  Finally, $P = P_{[1, \infty]}$ denotes the unique probability measure on $Z^\N$ with marginal $P_{[1,n]}$ on the first $n$ coordinates, which is guaranteed to exist by the Ionescu-Tulcea theorem.

%
%

In anticipation of our main abstract theorems, we define two technical assumptions about the Markov transition kernel $P_1$. The first pertains to contraction in negative regularity by averaging over compositions with random diffeomorphisms determined by time-1 evolution of the Markov chain.

\begin{assumption}[Time-1 Negative-Regularity Contraction]
\label{ass:negative_regularity_mixing}
Let $\{\varphi_z\}_{z\in Z} \subseteq \mathrm{Diff}^1(M)$ be a family of diffeomorphisms such that $z \mapsto \varphi_z$ is measurable.
We say that $P_1$ satisfies time-1 negative regularity contraction if there exist $\delta,\beta \in (0, 1)$, and $p\in (0,\infty)$ and a measurable mapping $V:Z\to [0,\infty)$ such that, for any $z \in Z$ and $\theta \in H^{-\delta}_0(M)$,
  \begin{equation*}
    \Big(\int_Z \|(\varphi_w)_* \theta\|_{H^{-\delta}(M)}^p\, \dif P_1(z)(w)\Big)^{1/p} \leq \beta e^{V(z)}\|\theta\|_{H^{-\delta}(M)}.
  \end{equation*}
\end{assumption}

In order to prove negative regularity mixing, we also need to assume that the Lyapunov function $V$ satisfies the following compatibility condition.

\begin{assumption}[Compatible Super-Lyapunov Function]
\label{ass:super_lyapunov}
Suppose that $P_1$ satisfies
Assumption~\ref{ass:negative_regularity_mixing} above for some Lyapunov function $V$, a family $\{\varphi_z\}$, and constants $\delta,\beta\in (0,1)$ and $p\in (0,\infty)$. We say that $V$ is a compatible super-Lyapunov function if, for some constants $K, \alpha > 1$ such that $
\beta K^{1/(\alpha-1)}<1
$,
  \begin{equation}
  \label{eq:superLyapunov}
    \int_Z e^{\alpha V(w)}\,\dif P_1(z)(w) \leq Ke^{V(z)} \quad \text{for any $z \in Z$.}
  \end{equation}
\end{assumption}

\begin{remark}\label{rem:neg-reg-equivalence}
We note that the abstract setting of~\cite{BFPS24} is more restrictive than the setting here. Indeed, while~\cite{BFPS24} appears to allow for more general random dynamical systems, since the resulting system must be Markovian anyway, we can handle this setting up to appropriately defining  $Z$, $P_1$ and $\{\varphi_z\}_{z\in Z}$ as they appear in this paper.

Recast in the notation of this paper, the conclusion of the main theorem of~\cite{BFPS24} implies that for some $V$ satisfying~\eqref{eq:superLyapunov} for some $K$, for each $p \in (0,2)$, there exists $C_p,\mu >0$ such that
\begin{align}\label{eq:BFPSmixing}
\left(\int_{Z^n} \|(\varphi_{z_n})_* \dots (\varphi_{z_1})_* \theta\|_{H^{-\delta}(M)}^p\, \dif P_{[1,n]}(z_0)(z_1, z_2, \dots, z_n)\right)^\frac{1}{p} \leq C_p e^{V(z_0) - \mu n} \|\theta\|_{H^{-\delta}(M)},
\end{align}
for all $n \in \mathbb N$. This estimate alone does not imply Assumptions~\ref{ass:negative_regularity_mixing} and~\ref{ass:super_lyapunov}. However, taking $N \in \mathbb N$ sufficiently large, we can verify these assumptions for the time-$N$ transition kernel, viewed on the state space $Z^N$. Indeed, Assumption~\ref{ass:negative_regularity_mixing} can be made true for any $p \in (0,2)$ and $\beta \in(0,1)$ by taking $N$ large enough for $C_pe^{-\mu N} < \beta$. By iterating~\eqref{eq:superLyapunov}, one observes that
\[
\mathbf E[e^{\alpha V(z_N)}] \leq K^{\frac{\alpha}{\alpha-1}}e^{V(z_0)}.
\]
Thus, for the time-$N$ Markov chain, Assumption~\ref{ass:super_lyapunov} holds with $K$ replaced by $K^{\frac{\alpha}{\alpha-1}}$, and $\beta = \beta(N)$ chosen sufficiently small.

On the other hand, we show in Theorem~\ref{thm:shortTimeToLong} below that Assumptions ~\ref{ass:negative_regularity_mixing} and ~\ref{ass:super_lyapunov} are sufficient to prove~\eqref{eq:BFPSmixing}, albeit with a smaller exponent $q < p$. Thus, Assumptions~\ref{ass:negative_regularity_mixing} and~\ref{ass:super_lyapunov} are essentially equivalent to~\eqref{eq:BFPSmixing}.
However, unlike~\eqref{eq:BFPSmixing}, it is easier to show that Assumptions~\ref{ass:negative_regularity_mixing} and~\ref{ass:super_lyapunov} are stable under perturbations of $P_1$ and $\{\varphi_z\}_{z\in Z}$.
\end{remark}

\subsection{Statement of the abstract results}\label{s:abstract_theory}
In this subsection, we state the main result in the abstract setting.
To formulate our first result, we introduce a Wasserstein-type functional that captures possible perturbations of the kernel and of flows. Let $m\geq 0$. We denote by $\Wasl_{d\wedge 1,m}$ the following weighted and truncated Wasserstein functional,
\begin{equation}
  \label{eq:wassersteinDef}
  \Wasl_{d\wedge 1,m}[\mu,\nu] \coloneq \inf_{\pi \in \Gamma(\mu, \nu)} \int_{\Diff^1(M)^2} (1+ \BiLip(\varphi)+ \BiLip(\psi))^m(\dist^1 (\varphi,\psi)\wedge 1)\, \dif \pi(\varphi,\psi)
\end{equation}
where $\mu, \nu \in \Prob(\mathrm{Diff}^1(M))$ and
\[
\Gamma(\mu, \nu) \coloneq\{\pi \in \Prob(\mathrm{Diff}^1(M)^2)\,:\, \pi(\cdot\times \mathrm{Diff}^1(M))=\mu ,\  \pi (\mathrm{Diff}^1(M)\times \cdot)=\nu\}\]
denotes the set of $(\mu, \nu)$-couplings. Recall that $\dist^1$ is the distance introduced in~\eqref{eq:diffDistance}.



In the first theorem of this section, we show that control on the Wasserstein functional defined above implies finite time stability of decay.

\begin{theorem}[Stability of unit-time decay]\label{thm:mainStability}
Let $P_1: Z\to \Prob (Z)$ be a kernel satisfying Assumption~\ref{ass:negative_regularity_mixing} for some $V$, $\delta,\beta\in (0,1)$ and $p\geq 1$, and a measurable family $(\varphi_z )_{z\in Z}\subseteq \Diff^1(M)$. Then, set $m= p(\frac{d}{2}+\delta)+1$.

Let $Q_1:Z \to \Prob(Z)$ be another measurable kernel and $(\psi_z )_{z\in Z}\subseteq \Diff^1(M) $ another measurable family of diffeomorphisms, not necessarily satisfying Assumption~\ref{ass:negative_regularity_mixing}.
For each $z \in Z$, define the measures $\mu_z,\nu_z \in \Prob(\mathrm{Diff}^1(M))$ by
\[
\mu_z(A) = \int 1_{A}(\varphi_w) \dif P_1(z)(w), \qquad \nu_z(A) = \int 1_{A}(\psi_w) \dif Q_1(z)(w)
\]
Then, we assume that for some $\varepsilon> 0$, for all $z \in Z$,
\begin{equation}
\label{eq:smallness_condition_Wasmodified}
 \Wasl_{d\wedge 1,m}[\mu_z,\nu_z] \leq \varepsilon^p e^{pV(z)}.
\end{equation}
Then the perturbed system $(Q_1, (\psi_z)_{z\in Z})$ satisfies the following negative-regularity decay estimate:
  \begin{equation}
  \label{eq:mainStability}
    \left(\int_Z \|(\psi_w)_* \theta\|_{H^{-\delta}(M)}^p\, \dif Q_1(z)(w)\right)^{1/p} \leq \beta' e^{V(z)}\|\theta\|_{H^{-\delta}(M)},
  \end{equation}
  for any $z \in Z$ and $\theta \in H^{-\delta}_0(M)$,
  where $\beta'= \beta+ C_1\varepsilon$ and $C_1$ depends only on $M,p$ and $\delta$.
\end{theorem}

\begin{remark}[Verifying the smallness condition~\eqref{eq:smallness_condition_Wasmodified}]
\label{r:Was_condition_checked}
In the setting of Theorem~\ref{thm:mainStability}, for any coupling $\Pi_z\in\Gamma(P_1(z),Q_1(z))$, the
pushforward measure $\widehat{\Pi}_z \coloneq (\varphi,\psi)_{\#}\Pi_z$ is a coupling of $\mu_z$ and $\nu_z$, and therefore
\begin{align}
\label{eq:condition_smallness_reduced}
\Wasl_{d\wedge 1,m}[\mu_z,\nu_z]
\leq \int_{Z^2} (1+ \BiLip(\varphi_{z_0})+ \BiLip(\psi_{z_1}))^m(\dist^1 (\varphi_{z_0},\psi_{z_1})\wedge 1)\, \dif \Pi_z(z_0,z_1).
\end{align}
In the later applications, we use a natural coupling $\Pi_z$ obtained by taking the same realization of the underlying randomness.
For example in the application to enhanced dissipation (Section~\ref{sec:enhanced-dissipation}), we use $P_1=Q_1$ and the corresponding diagonal coupling $\Pi_z=(\mathrm{Id},\mathrm{Id})_{\#}P_1(z)$. The above condition simplifies to

\begin{equation*}
\Wasl_{d\wedge 1,m}[\mu_z,\nu_z]
\leq \int_{Z} (1+ \BiLip(\varphi_{z'})+ \BiLip(\psi_{z'}))^m(\dist^1 (\varphi_{z'},\psi_{z'})\wedge 1)\, \dif P_1(z)(z').
\end{equation*}
In Sections~\ref{s:thermal_flows} and~\ref{sec:regularized}, the transition kernels differ (i.e., $P_1\neq Q_1$),
but the map assigning a Lagrangian flow to a state is the same for both systems ($\varphi=\psi$). In these cases, we use~\eqref{eq:condition_smallness_reduced} with the coupling $\Pi_z$ induced by solutions of the two stochastic systems driven by the same noise (a synchronous coupling). Thus, verifying~\eqref{eq:smallness_condition_Wasmodified} reduces to controlling the weighted $C^1$-distance between the coupled Lagrangian flows (see Lemma~\ref{lem:velocityToFlowMap} for the case of diffeomorphisms induced by ODEs).
\end{remark}

In the second theorem of this section, we give sufficient criteria for extending the finite time decay afforded by the previous theorem to arbitrary times.

\begin{theorem}[From unit-time to long-time decay]\label{thm:shortTimeToLong}
Let $(P_n)_{n\in\N}$ and $(\varphi_z)_{z\in Z}$ be a Markov kernel and a family of diffeomorphisms, respectively. Suppose that $P_1$ satisfies
Assumptions~\ref{ass:negative_regularity_mixing} and~\ref{ass:super_lyapunov} for some $V : Z \to [0,\infty)$, $\delta,\beta\in (0,1)$, $K > 1$, $p\in (0,\infty)$, a measurable family $\{\varphi_z\}_{z\in Z}$.

Let $q\in (0,p)$ satisfy
\begin{equation}
\label{eq:condition_pq}
\frac{pq}{p-q}\leq \alpha-1 \quad \Leftrightarrow \quad q \leq \frac{p(\alpha -1)}{p + \alpha - 1}.
\end{equation}
Then there exist $\gamma,\chi, C>0$ depending only on $M,\alpha,K,\beta,V,\delta,p$ and $q$ for which the following long-time decay estimate holds for any $z_0 \in Z$ and $\theta \in H^{-\delta}_0(M)$:
  \begin{align}
  \label{eq:shortTimeToLong}
    \left(\int_{Z^n} \|(\varphi_{z_n})_* \dots (\varphi_{z_1})_* \theta\|_{H^{-\delta}(M)}^q\, \dif P_{[1,n]}(z_0)(z_1, z_2, \dots, z_n)\right)^{1/q}&\\
    \nonumber
    \leq  C e^{-\gamma n} \exp(\chi V(z_0)) \|\theta\|_{H^{-\delta}(M)}&.
  \end{align}
Consequently, there exists a random constant $D = D(\theta, (z_n)_{n\in \N}) > 0$ satisfying $\int_{Z^\N} D^q \dif P_{[1,\infty]} \leq C$ such that
    \begin{equation}
    \label{eq:longNegRegMixingAS}
    \|(\varphi_{z_n})_*\dots (\varphi_{z_1})_* \theta\|_{H^{-\delta}(M)}
    \leq e^{-\gamma n}D\exp(\chi V(z_0))\|\theta\|_{H^{-\delta}(M)}.
  \end{equation}
  Furthermore, there is a constant $k > 0$ and a random constant $\tilde{D} = \tilde{D}((z_n)_{n \in \N}) > 0$, independent of the initial data $\theta$, satisfying $\int_{Z^\N} \tilde{D}^q \dif P_{[1,\infty]} \leq C$ and
    \begin{equation}
    \label{eq:longPosRegMixingAS}
    \|(\varphi_{z_n})_*\dots (\varphi_{z_1})_* \theta\|_{H^{-\delta}(M)}
    \leq e^{-\gamma n}\tilde{D}\exp(\chi V(z_0))\|\theta\|_{H^{k}(M)}.
  \end{equation}
\end{theorem}
\smallskip

\begin{remark}
  If the family $(\varphi_z)_z$ preserves volume, then $k > 0$ in~\eqref{eq:longPosRegMixingAS} can be taken to be any positive number by interpolation with $L^2(M)$, at the cost of choosing $\gamma > 0$ even smaller.
\end{remark}
\begin{remark}
While Theorems~\ref{thm:mainStability} and~\ref{thm:shortTimeToLong} are stated as separate results for technical reasons, they are applied jointly in every application presented in this paper.
\end{remark}

\section{Proof of Theorems~\ref{thm:mainStability} and~\ref{thm:shortTimeToLong}}
The proofs of Theorems~\ref{thm:mainStability} and~\ref{thm:shortTimeToLong} are given in Subsections~\ref{ss:proof_mainStability} and~\ref{ss:proof_shortTimeToLong}, respectively.

\subsection{Proof of Theorem~\ref{thm:mainStability}}
\label{ss:proof_mainStability}
To prove the stability result of Theorem~\ref{thm:mainStability}, we need the following result on composition of Sobolev distributions.

\begin{lemma}[Estimates for compositions of Sobolev distributions]
\label{lem:HDeltaPerturb}
There exists $C>0$ such that, for any diffeomorphisms $\varphi, \psi \in \Diff(M)$, any $\delta\in [0,1]$, and any $\theta \in H^{-\delta}(M)$, it holds that
  \begin{align}
  \label{eq:H-delta-interp}
    \|\varphi_* \theta\|_{H^{-\delta}(M)} &\leq \|\theta\|_{H^{-\delta}}\BiLip(\varphi)^{d/2 + \delta},\\
    \label{eq:HDeltaPerturb}
    \|\psi_* \theta\|_{H^{-\delta}(M)} &\leq \left(1 + C\left(\BiLip(\psi) + \BiLip(\varphi)\right)d^1(\varphi, \psi)\right)^{d/2+\delta}\|\varphi_* \theta\|_{H^{-\delta}(M)}.
  \end{align}
\end{lemma}

As will become clear from the proof of Theorem~\ref{thm:mainStability} below, the fact that~\eqref{eq:HDeltaPerturb} holds with unit constant when $\varphi=\psi$ plays a fundamental role in our perturbation argument. This property follows rather directly from the choice of the norm~\eqref{eq:definition_fractional_Sob_spaces} for the negative fractional Sobolev spaces.

\begin{proof}
We only prove~\eqref{eq:HDeltaPerturb}, as~\eqref{eq:H-delta-interp} is similar and simpler.
Set $h\coloneq \psi\circ\varphi^{-1}$, and note that $\psi_* \theta= h_* (\varphi_* \theta)$. As $h\in \Diff(M)$ is onto, it follows that $\operatorname{Lip}(h)\geq 1$.
To prove the result, we use complex interpolation of linear operators~\cite[Theorem C.2.6]{Analysis1} on the interpolation spaces~\eqref{eq:definition_fractional_Sob_spaces}. Thus, it suffices to show that the linear operator $h_*$  satisfies the bounds
  \begin{equation}
    \label{eq:HDeltaPerturb_claim}
    \|h_* \theta\|_{H^{-k}(M)} \leq \left(1 + C\left(\BiLip(\psi) + \BiLip(\varphi)\right)d^1(\varphi, \psi)\right)^{d/2+k}\| \theta\|_{H^{-k}(M)},
  \end{equation}
  for $k\in \{0,1\}$.
To prove~\eqref{eq:HDeltaPerturb_claim}, first note that, by a change of variables,
$$
\|h_*\theta\|_{L^2(M)}
\leq \|J_{h^{-1}}\|_{L^\infty(M)}^{1/2}\|\theta\|_{L^2(M)}
\leq \operatorname{Lip}(h^{-1})^{d/2}\|\theta\|_{L^2(M)}.
$$
Second, by duality and the chain rule,
\begin{align*}
\|h_*\theta\|_{H^{-1}(M)}
&\leq \|\theta\|_{H^{-1}(M)}
\sup_{\|f\|_{H^1(M)}\leq 1}\|f\circ h\|_{H^1(M)} \\
&\leq\|\theta\|_{H^{-1}(M)} \operatorname{Lip}(h)\,
\|J_{h^{-1}}\|_{L^\infty(M)}^{1/2} \\
&\leq \|\theta\|_{H^{-1}(M)} \operatorname{Lip}(h)\,
\operatorname{Lip}(h^{-1})^{d/2}.
\end{align*}
Moreover, by the definition of $\dist^1$ in~\eqref{eq:diffDistance},
$$
\operatorname{Lip}(h)
\leq 1+C\operatorname{BiLip}(\varphi)\dist^1(\varphi,\psi),
\qquad
\operatorname{Lip}(h^{-1})
\leq 1+C\operatorname{BiLip}(\psi)\dist^1(\varphi,\psi).
$$
This completes the proof of~\eqref{eq:HDeltaPerturb_claim} in the extreme cases $k\in \{0,1\}$.
\end{proof}

With the above result at our disposal, we are ready to prove Theorem~\ref{thm:mainStability}.

\begin{proof}[Proof of Theorem~\ref{thm:mainStability}]

Fix $z\in Z$ and set $\mu=\mu_z$, $\nu=\nu_z$. Let $s\coloneq\frac{d}{2}+\delta$. For $\xi\in \mathrm{Diff}^1(M)$, we set
$$
N_\xi\coloneq\|(\xi)_*\theta\|_{H^{-\delta}(M)},
\qquad
B_\xi\coloneq\BiLip(\xi).
$$
By slight abuse of notation, we use $(\varphi,\psi)$ to denote a sample from the coupling $\pi \in \Gamma(\mu,\nu)$.
We claim that there exists $C\geq 1$ independent of $\varphi,\psi$ such that for all $\varphi,\psi\in\Diff^1(M)$,
\begin{equation}
\label{eq:claim_N_distance_claim}
|N_\varphi^p-N_\psi^p|
\leq C\|\theta\|_{H^{-\delta}(M)}^p
(B_\varphi+B_\psi)^{ps+1}\,
(\dist^1(\varphi,\psi)\wedge 1).
\end{equation}
Let us first prove the above in the case $(B_\varphi+B_\psi)(\dist^1(\varphi,\psi)\wedge 1)\leq 1$.  Since $B_\varphi+B_\psi \geq 4$, it follows that $\dist^1(\varphi,\psi) = \dist^1(\varphi,\psi) \wedge 1 \leq \frac{1}{4}$.
Then, it follows from Lemma~\ref{lem:HDeltaPerturb} applied with $\varphi$ and $\psi$ in both orders that
\begin{align*}
|N_\varphi-N_\psi|
&\leq
\big[
\left(1+C(B_\varphi+B_\psi)\dist^1(\varphi,\psi)\right)^s-1
\big](N_\varphi+N_\psi)\\
&\leq
C(B_\varphi+B_\psi)(\dist^1(\varphi,\psi)\wedge 1)(N_\varphi+N_\psi),
\end{align*}
where we used the elementary inequality $(1+Cx)^s-1\leq ((C+1)^s -1)x$ for $x\leq1$.

Combining this with~\eqref{eq:H-delta-interp}, we get
\begin{align}
\label{eq:N_distance_p_bound}
|N_{\varphi}^p-N_{\psi}^p|
&\leq C (N_{\varphi}^{p-1}+N_{\psi}^{p-1})|N_{\varphi}-N_{\psi}|\\
\nonumber
&\leq C (B_\varphi+B_\psi)\,(\dist^1(\varphi,\psi)\wedge 1) (N_{\varphi}^{p}+N_{\psi}^{p})\\
\nonumber
&\leq C\|\theta\|_{H^{-\delta}(M)}^p
(B_\varphi+B_\psi)^{ps+1}\,
(\dist^1(\varphi,\psi)\wedge 1),
\end{align}
as desired.

Alternatively, if $(B_\varphi+B_\psi)(\dist^1(\varphi,\psi)\wedge 1)> 1$, we simply apply~\eqref{eq:H-delta-interp} and the triangle inequality to get
\begin{align*}
|N_{\varphi}^p-N_{\psi}^p| &\leq C \|\theta\|_{H^{-\delta}}^p(B_\varphi^{ps} + B_\psi^{ps})\\
&\leq  C\|\theta\|_{H^{-\delta}}^p(B_\varphi + B_\psi)^{ps +1} (\dist^1(\varphi,\psi) \wedge 1).
\end{align*}

Recalling that $\pi\in \Gamma(\mu,\nu)$, we have
\begin{align*}
&\Big(\int_{\Diff^1(M)} N_\psi^p\,\dif \nu(\psi)\Big)^{1/p}\\
&=
\Big(\int_{\Diff^1(M)^2} N_\psi^p\,\dif \pi (\varphi,\psi)\Big)^{1/p}\\
&=
\Big(\int_{\Diff^1(M)^2} \big(N_\varphi^p + N_\psi^p-N_\varphi^p\big)\,\dif \pi (\varphi,\psi)\Big)^{1/p} \\
&\leq
\Big(\int_{\Diff^1(M)} N_\varphi^p\,\dif \mu(\varphi)\Big)^{1/p} +
\Big(
\int_{\Diff^1(M)^2} |N_\varphi^p-N_\psi^p|\,\dif\pi(\varphi,\psi)
\Big)^{1/p}  \\
&\stackrel{(i)}{\leq}
\beta\exp(V(z))\|\theta\|_{H^{-\delta}(M)} \\
&\quad +C\|\theta\|_{H^{-\delta}(M)}
\Big(
\int_{\Diff^1(M)^2}
(B_\varphi+B_\psi)^{ps+1}\,
(\dist^1(\varphi,\psi)\wedge 1)\,\dif\pi (\varphi,\psi)
\Big)^{1/p}
\end{align*}
where in $(i)$ we used Assumption~\ref{ass:negative_regularity_mixing} and the bound~\eqref{eq:claim_N_distance_claim}.
Taking the infimum over all such $\pi \in \Gamma(\mu,\nu)$ yields
$$
\Big(\int_{\Diff^1(M)} N_\psi^p\,\dif \nu(\psi)\Big)^{1/p}
\leq\|\theta\|_{H^{-\delta}(M)}
(\beta\exp(V(z)) +C   \Wasl_{d\wedge 1,m}[\mu,\nu]^{1/p}).
$$
The conclusion now follows by using~\eqref{eq:smallness_condition_Wasmodified}.
\end{proof}

\subsection{Proof of Theorem~\ref{thm:shortTimeToLong}}
\label{ss:proof_shortTimeToLong}
The proof of Theorem~\ref{thm:shortTimeToLong} is based on an induction argument together with the Markov property of $(P_n)_{n\in \N}$.

\begin{proof}[Proof of Theorem~\ref{thm:shortTimeToLong}]
  For notational convenience, we let
  $$
  \chi\coloneq\frac{\alpha}{\alpha-1}.
  $$
  Moreover, from the assumption $\beta K^{1/(\alpha-1)}<1$, there exists $\g>0$ such that
  \begin{align}\label{eq:gammagap}
  \beta K^{1/(\alpha-1)} e^{2\g}<1.
\end{align}
  We prove~\eqref{eq:shortTimeToLong} with $C=1$ and decay rate $2\gamma$ by induction.
  The case $n=1$ of the latter follows from Assumption~\ref{ass:negative_regularity_mixing} as $\beta e^{V(z_0)}\leq e^{-2\g } e^{\chi V(z_0)}$ as $K\geq 1$ and the above choice of $\chi$ and $\g$.
  Next, assume that this stronger estimate holds for some $n\geq 1$.
 Using the Markov property, we apply the inductive step with $\theta$ replaced by $(\varphi_{z_1})_* \theta$, which also has mean zero. Combining this with H\"older's inequality, we obtain
  \begin{align*}
    &\Big(\int_{Z^{n+1}} \|(\varphi_{z_{n+1}})_* \dots (\varphi_{z_1})_* \theta\|_{H^{-\delta}(M)}^q \dif P_{[1,n+1]}(z_0)(z_1, \dots, z_{n+1})\Big)^{1/q}\\
    &\qquad  = \Big[\int_Z  \Big( \int_{Z^n}  \|(\varphi_{z_{n+1}})_* \dots (\varphi_{z_1})_* \theta\|_{H^{-\delta}(M)}^q \dif P_{[1,n]}(z_1)(z_2, \dots, z_{n+1}) \Big)\,\dif P_1(z_0)(z_1) \Big]^{1/q}\\
    &\qquad \leq e^{-2\gamma n} \Big( \int_Z e^{q\chi V(z_1)} \|(\varphi_{z_1})_* \theta\|_{H^{-\delta}(M)}^q \,\dif P_1(z_0)(z_1)\Big)^{1/q}\\
    &\qquad \leq e^{-2\gamma n}\Big(\int_Z e^{\lambda  V(z_1)}\,\dif P_1(z_0)(z_1)\Big)^{\frac{p-q}{pq}}	\Big(\int_Z \|(\varphi_{z_1})_* \theta\|_{H^{-\delta}(M)}^{p}\, \dif P_1(z_0)(z_1) \Big)^{1/p}.
  \end{align*}
  Above, we set $\lambda\coloneq\frac{\chi pq}{p-q} = \frac{\alpha}{\alpha-1} \cdot \frac{pq}{p-q}$.
 By~\eqref{eq:condition_pq}, we have $\lambda \leq \alpha$. Hence, it follows from~\eqref{eq:superLyapunov} and Jensen's inequality that
  \begin{align*}
  \int_Z e^{\lambda  V(z_1)}\,\dif P_1(z_0)(z_1)&\leq\Big(
  \int_Z e^{\alpha  V(z_1)}\,\dif P_1(z_0)(z_1) \Big)^{\frac{\lambda}{\alpha}} \\
  &\leq  \left(K\exp(V(z_0))\right)^{\frac{\lambda}{\alpha}} = \left(K\exp(V(z_0))\right)^{\frac{pq}{(\alpha -1)(p-q)}}.
  \end{align*}
Equivalently,
\begin{align*}
\left(\int_Z e^{\lambda  V(z_1)}\,\dif P_1(z_0)(z_1)\right)^{\frac{p-q}{pq}} \leq \left(Ke^{V(z_0)}\right)^{\frac{1}{\alpha -1}}.
\end{align*}
Combining these estimates, and again using Assumption~\ref{ass:negative_regularity_mixing}, we obtain
\begin{align*}
&\left(\int_{Z^{n+1}} \|(\varphi_{z_{n+1}})_* \dots (\varphi_{z_1})_* \theta\|_{H^{-\delta}(M)}^q \dif P_{[1,n+1]}(z_0)(z_1, \dots, z_{n+1})\right)^{1/q}\\
&\quad\leq e^{-2\g n} \beta K^{\frac{1}{\alpha-1}} e^{\left(1+\frac{1}{\alpha-1}\right) V(z_0)} \|\theta\|_{H^{-\delta}(M)}\\
&\quad = e^{-2\gamma n } \beta K^{\frac{1}{\alpha-1}} e^{\chi V(z_0)} \|\theta\|_{H^{-\delta}(M)}
\end{align*}
Applying~\eqref{eq:gammagap}, we can bound the above by $e^{-2\gamma (n+1)}e^{\chi V(z_0)}\|\theta\|_{H^{-\delta}(M)}$, showing that
the stronger estimate holds with $n$ replaced by $n+1$. This completes the inductive step.

We now prove the pointwise bound~\eqref{eq:longNegRegMixingAS}, keeping the above choices of $\gamma,\chi>0$.
The result is immediate in the case $\theta = 0$, so we assume $\|\theta\|_{H^{-\delta}} > 0$.
  Define
  \begin{equation*}
    D_n \coloneq e^{\gamma n}e^{-\chi V(z_0)}\|\theta\|_{H^{-\delta}(M)}^{-1}\|(\varphi_{z_n})_* \dots (\varphi_{z_1})_* \theta\|_{H^{-\delta}(M)}.
  \end{equation*}
  Then the stronger estimate implies
  \begin{align*}
\int D^q_n \dif P_{[1, \infty]} (z_0)
&\leq e^{q\gamma n}\exp(-q\chi V(z_0))\|\theta\|_{H^{-\delta}(M)}^{-q}\\
&\qquad \times \int_{Z^n}\|(\varphi_{z_n})_* \dots (\varphi_{z_1})_* \theta\|_{H^{-\delta}(M)}^q\dif P_{[1,n]}(z_0)(z_1, \dots, z_n) \leq e^{-q n\g}.
\end{align*}
Then~\eqref{eq:longNegRegMixingAS} follows from the above by letting $D \coloneq \sup_{n \in \N} D_n$ and using $D^q\leq \sum_{n\in \N} D_n^q$.

We now deduce~\eqref{eq:longPosRegMixingAS} from~\eqref{eq:longNegRegMixingAS}.
Let $\{e_n\}_{n \in \N}$ be an orthonormal basis of eigenvectors of the Laplacian for $L^2_0(M)$, ordered by increasing eigenvalue.
Let $D_n$ be the constant satisfying~\eqref{eq:longNegRegMixingAS} for initial data $\theta = e_n$.
By Weyl's law~\cite[Theorem~1.1]{Hormander1968}, we have $\|e_n\|_{H^k} \approx n^{k/d}$ for $k \geq 1$.
Let $\theta = \sum_n \alpha_n e_n \in H^k(M)$.
Then
\begin{equation*}
  \sum_n \alpha_n^2 \|e_n\|_{H^k}^2 \approx \sum_n \alpha_n^2 n^{2k/d} < \infty.
\end{equation*}
By~\eqref{eq:longNegRegMixingAS} we have
\begin{align*}
  \|(\varphi_{z_n})_*\dots (\varphi_{z_1})_* \theta\|_{H^{-\delta}}
  &\leq e^{-\gamma n}\exp(\chi V(z_0)) \sum_i D_i |\alpha_i| \|e_i\|_{H^{-\delta}}\\
  &\leq e^{-\gamma n}\exp(\chi V(z_0))\|\theta\|_{H^k} C\sum_i D_i i^{-k/d}.
\end{align*}
Set $\tilde{D} \coloneq C\sum_n D_n n^{-k/d}$ and choose $k > d/\min\{1,q\}$ to conclude.
\end{proof}

\section{Enhanced dissipation}\label{sec:enhanced-dissipation}


Fix $U = C^2_\text{sol}(\T^d; \R^d)$ to be the space of $C^2$ divergence-free vector fields on $\mathbb T^d$, and let $(u_t)_{t\geq 0} \in C([0,\infty);U)$ be a continuous time Markov process on $U$. Letting $\kappa \in[0,1]$, we consider solutions $\theta^\kappa_t$ to the following advection-diffusion equation
\begin{equation}
  \label{eq:advectionDiffusion}
  \begin{cases}
    \partial_t \theta^\kappa_t + u_t \cdot \nabla \theta^\kappa_t = \kappa \Delta \theta^\kappa_t, &\qquad \text{for $t > 0$,}\\
    \theta^\kappa_0 = \theta_0.
  \end{cases}
\end{equation}
Our main theorem is that if the above exhibits negative regularity mixing in the pure advective case $\kappa = 0$, then we have an enhanced dissipation estimate for $\theta^\kappa$ in $L^2$.

\begin{theorem}\label{thm:enhanced-dissipation}
Fix $p,\alpha >1$. Let $(u_t)_{t\geq 0}$ be a continuous time Markov process in $U$ (with $u_0 \in U$ fixed). Assume the existence of a function $V : U \to [0,\infty)$ satisfying a super-Lyapunov property
\[
\mathbf E[\exp(\alpha V(u_1))] \leq K \exp(V(u_0)).
\]
for some constant $ K \geq  1$.
Furthermore, assume that for every $N >0$, there exists $C_N >0$ such that
\begin{equation}\label{eq:exp-moment-base-velocity-new}
  \mathbf E\Big[ \exp\Big(N \int_0^N \|u_s\|_{C^{2}} \dif s\Big)\Big]  \leq C_N\exp(V(u_0)/2).
\end{equation}

Assume that there exist $\delta\in(0,1)$ and $C_0,\beta>0$ such that, for every $u_0\in U$ and $\theta_0\in H^{-\delta}_0(\mathbb T^d)$, the solution to~\eqref{eq:advectionDiffusion} with $\kappa = 0$ obeys the negative-regularity mixing estimate
\begin{align}\label{eq:cont_time_negative_regularity_mixing}
\mathbf E[\|\theta^0_t\|_{H^{-\delta}}^p]^{\frac{1}{p}} \leq C_0 e^{-\beta t}e^{V(u_0)}\|\theta_0\|_{H^{-\delta}}, \qquad t\geq 0.
\end{align}
 Then, there exist $\gamma >0$ and $\chi \geq 1$ such that for all $\kappa >0$ sufficiently small and initial conditions $\theta_0 \in L^2_0(\mathbb T^d)$, solutions to~\eqref{eq:advectionDiffusion} exhibit the enhanced dissipation estimate
\begin{align}
\|\theta^\kappa_t\|_{L^2} \leq D\kappa^{-{\delta/2}} e^{-\gamma t}e^{\chi V(u_0)}\|\theta_0\|_{L^2}, \qquad t\geq 0.
\end{align}
Above, $D = D(\kappa, (u)_{t\geq 0})$ is a random constant with $\E D \leq C$, where $C > 0$ is an absolute constant.
\end{theorem}

In the following we collect some examples of continuous Markov processes $(u_t)_{t\geq 0}$ to which the above result applies. Although we assume for simplicity that the velocity $(u_t)_{t \geq 0}$ is Markov in the statement of the theorem, the proof only requires the weaker assumption that unit-time chunks of the velocity form a Markov process

\begin{example}\
\begin{itemize}
    \item {\rm (Stochastic fluids)}
Following Remark~\ref{rem:neg-reg-equivalence}, Theorem~\ref{thm:enhanced-dissipation} applies to any process $(u_t)_{t\geq 0}$, with enough regularity to satisfy~\eqref{eq:exp-moment-base-velocity-new}, for which negative regularity mixing in the sense of~\cite{BFPS24} is known. In particular, one can take $(u_t)_{t\geq 0}$ to be a solution to the 2D Navier-Stokes equation with nondegenerate stochastic forcing.
    \item {\rm (Composition of iid diffeomorphisms)}
     Suppose that the unit-time processes of $(u_t)_{t\geq 0}$
\[
(u_t)_{0\leq t < 1}, (u_t)_{1\leq t < 2},
(u_t)_{2\leq t < 3}, \ldots \]
are iid and satisfy~\eqref{eq:exp-moment-base-velocity-new}. Hence, the time-$n$ flow map of $(u_t)_{t\geq 0}$ is a composition of $n$ iid diffeomorphisms for any $n = 1,2,3,\ldots$. If the time-$n$ flow maps satisfy Assumptions 2 and 3 in~\cite{BFPS24} (one can take $V \equiv 0$ for such examples, so Assumption 1 is satisfied automatically), then the hypotheses of Theorem~\ref{thm:enhanced-dissipation} are satisfied.
\item {\rm (Pierrehumbert flows)}
One particular example of compositions of random maps exhibiting negative regularity mixing are the Pierrehumbert flows on $\mathbb T^2$. While this flow is usually discrete in time, we can construct the flow from a  continuous time velocity field $u_t$ as follows. Let $X_0,Y_0,X_1,Y_1,\ldots$ be a sequence of iid $\mathrm{Unif}([0,1])$ random variables, let $\tau > 0$ and let  $\chi \in C_c([0,1])$ be nonnegative and not identically zero. Then for each $t$, set
\[
u_t = \sum_{k=0}^\infty\left[ \chi(t - 2k) \begin{pmatrix} \tau\sin(2\pi (y - Y_k))\\
0
\end{pmatrix} + \chi(t-2k-1) \begin{pmatrix} 0\\
\tau \sin(2\pi (x - X_k))
\end{pmatrix}\right].
\]
Since all but one term vanishes on any interval $[n,n+1]$ for $n = 0,1,\ldots$, the time-$2n$ flow maps of $u_t$ are compositions of $n$ iid random diffeomorphisms.
At discrete times, this reproduces the same dynamics as the system considered, for instance, in~\cite{BlumenthalCotiZelatiGvalani2023}. Assumptions 2 and 3 of~\cite{BFPS24} were proved  for  Pierrehumbert flows in~\cite[Lemmas 4.1, 4.2 and 4.4 and Proposition 4.5]{BlumenthalCotiZelatiGvalani2023}.
\end{itemize}
\end{example}

\begin{proof}[Proof of Theorem~\ref{thm:enhanced-dissipation}] We break the proof up into a number of steps.

\textbf{Step 1} (reformulation in terms of unit time chunks): Let $S_t :U \to \mathcal P(U)$ be the time-$t$ Markov transition kernel associated with the process $(u_t)_{t\geq 0}$. For $u\in U$, write $\Law_u$ for the law of the process started at $u$, and set

\[
S_{[0,1]}(u) \coloneq \Law_u((u_s)_{0\leq s\leq1})\in \mathcal P(C([0,1];U)).
\]
Then, defining $Z \coloneq C([0,1];U) \times \{y\in C([0,1];\mathbb R^d):y_0=0\}$, we define a new transition kernel on $Z$ by
\[
P_1((u_t,y_t)_{t \in [0,1]}) = S_{[0,1]}(u_1) \otimes \mathbb W
\]
where $\mathbb W \in \mathcal P(C([0,1];\mathbb R^d))$ denotes the Wiener measure on Brownian paths starting from $0$.
 To distinguish between the state space $Z$ and the underlying process $(u_t,W_t)_{t\geq 0}$, we shall denote a generic unit-time sample path $z = (v_t,{y_t})_{t \in [0,1]}$.

 We then define the family of diffeomorphisms $\{\varphi^\kappa_z\}_{z\in Z}$. Given $z$, we solve the SDE for each initial condition $x \in \mathbb T^d$,
\begin{align}\label{eq:small_noise_sde}
 X_t^\kappa = x +\int_0^t v_{s}(X_s^\kappa) \dif s + \sqrt{2\kappa } y_t, \quad t \in [0,1],
\end{align}
Then, we set $\varphi_z^\kappa(x) \coloneq X_1^\kappa(x)$. By the stochastic method of characteristics, integrating over the Wiener measure recovers $\theta^\kappa_1$, taking $v = (v_t)_{t\in[0,1]} = (u_t)_{t\in [0,1]}$:
\begin{align}\label{eq:partial_averaging}
\theta_1^\kappa =  \int_{C([0,1];\mathbb R^d)} (\varphi_{(v,y)}^\kappa)_* \theta_0 \dif \mathbb W(y).
\end{align}
Iterating the above for compositions $(\varphi_{z_n})_*\cdots(\varphi_{z_1})_*\theta_0$, we get similar formulas for each $\theta_n^\kappa$.
Note that for the $\kappa = 0$ case, $\varphi_{(v,y)}$ is constant with respect to $y$.

Next, by slight abuse of notation, we extend $V$ to $Z$ by
\[
V((v_t,y_t)_{t\in[0,1]}) \coloneq V(v_1).
\]
Then, in the case $\kappa = 0$, by~\eqref{eq:cont_time_negative_regularity_mixing}, without loss of generality, we may re-scale $v_t$ so that there exist constants $\beta \in (0,1)$, $\alpha > 1$ and $K \geq 1$ such that $P_1, \{\varphi^{0}_z\}$ and $V$ satisfy Assumptions~\ref{ass:negative_regularity_mixing} and~\ref{ass:super_lyapunov} (see also Remark~\ref{rem:neg-reg-equivalence}).

\textbf{Step 2} (stability of decay for $\kappa > 0$):  with $P_1 = Q_1$,
and comparing the two families of diffeomorphisms $\{\varphi^0_z\}$ and $\{\varphi^\kappa_z\}$ for some $\kappa >0$ sufficiently small. We already know $(P_1,\{\varphi^0_z\})$ satisfies Assumptions~\ref{ass:negative_regularity_mixing} and~\ref{ass:super_lyapunov}.
We shall apply Theorem~\ref{thm:mainStability} to show that the same is true for $(P_1,\{\varphi^\kappa_z\})$.

 Indeed, fix $z = (v_t,y_t)_{t\in[0,1]} \in Z$ and define for each $\kappa \geq 0$
 \[
  \mu^\kappa(A)  = \int 1_A(\varphi_{w}^\kappa) \dif P_1(z)(w).
 \]
We then  estimate $\Wasl_{d\wedge 1,m}[\mu^0,\mu^\kappa]$. To do this,
we choose the coupling $\pi \in \Gamma(\mu^0,\mu^\kappa)$ whose sample flow maps $(\psi,\tilde \psi)$ are defined by solving the SDE~\eqref{eq:small_noise_sde} with zero and positive $\kappa$ respectively, but  with identical instances of $(v_{t})_{t\in[0,1]}$.
Take such a sample $(\psi,\tilde \psi)$ from $\pi$, and  corresponding solutions $(X_t^0,X^\kappa_t)_{t\in[0,1]}$ to the SDE~\eqref{eq:small_noise_sde}. Let $\overline \psi$ be the time-1 flow map for the vector field $(v_t( \cdot + \sqrt{2\kappa} y_t))_{t \in [0,1]}$. Observe that $\overline \psi$ is related to $\tilde \psi$ by
\[
\tilde \psi(x) =  X^\kappa_1 =\overline \psi(x) + \sqrt{2\kappa} y_1.
\]
We compare the differences of the  vector fields underlying $\psi$ and $\overline \psi$ by
\[
\|v_t-v_t( \cdot + \sqrt{2\kappa} y_t)\|_{L^\infty} \leq \sqrt{2\kappa} |y_t| \|\nabla v_t\|_{L^\infty}.
\]
Using this together with
 Lemma~\ref{lem:velocityToFlowMap}  yields
\begin{align*}
d^1(\psi,\overline \psi) &\leq C \left(\int_0^1 \sqrt{\kappa} |y_t| \|v_t\|_{C^1}\dif t\right)^{\frac{1}{2}} \exp\left( C \int_0^1\|v_t\|_{C^2}\dif t \right) \\
&\leq C\kappa^{\frac{1}{4}} \sup_{t\in[0,1]} \sqrt{|y_t|} \exp\left(C\int_0^1 \|v_t\|_{C^2} \dif t\right).
\end{align*}
On the other hand,
\[
d^1(\psi,\tilde \psi) \leq d^1(\psi,\overline \psi) + d^1(\tilde\psi,\overline \psi).
\]
Since $\tilde \psi$ and $\overline \psi$ differ only by a translation $\sqrt{2\kappa} y_1$, we apply  Lemma~\ref{lem:velocityToFlowMap} once more to get
\[
d^1(\tilde\psi,\overline \psi) \leq C {\kappa}^\frac{1}{2} | y_1| \bigl(1+d^2(\overline \psi,I)\bigr)\leq C\kappa^\frac{1}{2} \sup_{t \in[0,1]} |y_t| \exp  \left(C\int_0^1 \|v_t\|_{C^2}\dif t\right).
\]
Hence,
\[
d^1(\psi,\tilde \psi) \wedge 1 \leq C \kappa^{\frac{1}{4}} \sup_{t \in[0,1]} (1+|y_t|) \exp  \left(C\int_0^1 \|v_t\|_{C^2}\dif t\right)
\]
Also by Lemma~\ref{lem:velocityToFlowMap},
\[
\BiLip(\psi) + \BiLip(\tilde \psi) \leq C  \exp  \left(C\int_0^1 \|v_t\|_{C^2}\dif t\right).
\]
Hence, for any $m > 0$, by~\eqref{eq:exp-moment-base-velocity-new} and straightforward Gaussian estimates, there exists $C_m$ such that
\[
\Wasl_{d\wedge 1,m}[\mu^0,\mu^\kappa] \leq C_m \kappa^{1/4}\exp(V(v_1)/2) = C_m \kappa^{1/4}\exp(V(z)/2).
\]
Taking $m = p (\frac{d}{2} + \delta) + 1$, and $\eps = (C_m \kappa^{\frac{1}{4}})^\frac{1}{p}$, we can take $\kappa$ sufficiently small so that~\eqref{eq:mainStability} holds for $(P_1,\{\varphi_z^\kappa\})$ with a contraction factor $\beta'>0$ satisfying $\beta'K^{1/(\alpha-1)}<1$. Consequently, we may apply Theorem~\ref{thm:shortTimeToLong} with
\[
0< q < \min\left\{1,\frac{p(\alpha-1)}{p+\alpha -1} \right\}.
\]
This yields, for some $\chi > 0$ and $\gamma > 0$,
\begin{align*}
    &\left(\int_{Z^n} \|(\varphi_{z_n}^\kappa)_* \dots (\varphi_{z_1}^\kappa)_* \theta_0\|_{H^{-\delta}(M)}^q\, \dif P_{[1,n]}(z_0)(z_1, z_2, \dots, z_n) \right)^{\frac{1}{q}}&\\
    &
    \qquad \leq e^{-\gamma n} \exp(\chi V(z_0)) \|\theta_0\|_{H^{-\delta}(M)}.
\end{align*}
Above, we take $z_0 = (u_0, 0)_{t\in[0,1]}$.

Then, applying~\eqref{eq:partial_averaging},
\begin{align*}
  \E\left[\|\theta^\kappa_n\|_{H^{-\delta}}^q\right] &\leq \|\theta_0\|_{L^2}^{q(1-q)}\left(\int_{Z^n} \|(\varphi_{z_n}^\kappa)_* \dots (\varphi_{z_1}^\kappa)_* \theta_0\|_{H^{-\delta}(M)}^q\, \dif P_{[1,n]}(z_0)(z_1, z_2, \dots, z_n) \right)^{q}\\
                                                     &\leq e^{-\gamma n} \exp(\chi V(z_0))\|\theta_0\|_{L^2}^q,
\end{align*}
and interpolating (using $\|\theta^\kappa_n\|_{H^{-\delta}} \leq \|\theta_0\|_{L^2}$ almost surely) we conclude
\begin{equation*}
  \E\left[\|\theta^\kappa_n\|_{H^{-\delta}}\right] \leq e^{-\gamma n} \exp(\chi V(z_0))\|\theta_0\|_{L^2},
\end{equation*}
making $\gamma$ smaller and $\chi$ larger as needed.
Arguing similarly as in the end of the proof of Theorem~\ref{thm:shortTimeToLong} (and choosing $\gamma$ even smaller), there exist random constants $D_0, \widetilde D_0 >0$ with finite expectation such that
\[
\sup_{n \in \mathbb N,\, n\geq 1}e^{\gamma n}\|\theta^\kappa_n\|_{H^{-\delta}} \leq D_0\exp(C V(u_0))
\|\theta^\kappa_1\|_{H^1} \leq \kappa^{-1/2} \widetilde D_0\exp(C V(u_0)) \|\theta_0\|_{L^2}.
\]

\textbf{Step 3} (from negative regularity mixing to enhanced dissipation):
  for general $t \geq 1$, we observe
  \begin{align*}
    \|\theta^\kappa_t\|_{H^{-\delta}} &\leq C\|\theta^\kappa_{\floor{t}}\|_{H^{-\delta}} \exp\left(2\int_{\floor{t}}^t \Lip(u_r) \dif r\right) \\
    &\leq C\kappa^{-1/2}\widetilde D_0\exp(C V(u_0)) e^{-\gamma t} \|\theta_0\|_{L^2}\exp\left(2\int_{\floor{t}}^t \Lip(u_r) \dif r\right).
  \end{align*}
  On the other hand, energy estimates for $\theta^\kappa_t$ give
  \begin{equation*}
    \|\theta^\kappa_{t+s}\|_{H^1} \leq C\|\theta^\kappa_t\|_{H^1}\exp\left(2\int_t^{t+s}\Lip(u_r) \dif r\right),
  \end{equation*}
  and also that for every $n \in \N$ there is some $t \in [n, n+1]$ with
  \begin{equation*}
    \|\theta^\kappa_{t}\|_{H^1} \leq C\kappa^{-1/2}\|\theta_0\|_{L^2}.
  \end{equation*}
  Putting these together, we conclude that for every $t \geq 2$,
  \begin{equation*}
    \|\theta^\kappa_t\|_{H^1} \leq C\kappa^{-1/2}\|\theta_0\|_{L^2}\exp\left(2\int_{\floor{t}-1}^{\floor{t}+1}\Lip(u_r) \dif r\right).
  \end{equation*}

  By~\eqref{eq:exp-moment-base-velocity-new}, we have the expectation bound
  \begin{equation*}
  \E\left[\exp\left(2\int_{\floor{t}-1}^{\floor{t}+1} \Lip(u_r) \dif r\right)\right] \leq C  \E\left[ \exp(V(u_{\floor{t}-1})/2)\right] \leq C \exp(V(u_0)/2).
  \end{equation*}
  for some $C > 0$.
  Define
  \begin{equation*}
    D_1 \coloneq \sup_{t \in \N,\, t\geq 1} t^{-2} \exp(-V(u_0)/2)\exp\left(2\int_{t-1}^{t+1} \Lip(u_r) \dif r\right).
  \end{equation*}
  Then by~\eqref{eq:exp-moment-base-velocity-new}, for any $r \geq 1$,
  \begin{align*}
    \mathbf E[D_1^r]  &\leq \sum_{t=1}^\infty \mathbf E\left[ t^{-2r}\exp(-rV(u_0)/2)\exp\left(2r\int_{t-1}^{t+1} \Lip(u_r) \dif r\right)\right]\\
             &\leq C_r\sum_{t=1}^\infty t^{-2r}\\
             &\leq C_r.
  \end{align*}
  By the definition of $D_1$, we have
  \begin{equation*}
    \|\theta^\kappa_t\|_{H^1} \leq C D_1 \kappa^{-1/2} t^2 \exp(V(u_0)/2) \|\theta_0\|_{L^2}
  \end{equation*}
  and
  \begin{equation*}
    \|\theta^\kappa_t\|_{H^{-\delta}} \leq C\kappa^{-1/2}\widetilde D_0D_1 t^2 \exp\left(\left(C + \frac{1}{2}\right)V(u_0)\right)e^{-\gamma t} \|\theta_0\|_{L^2}
  \end{equation*}
  Interpolating between $H^{-\delta}$ and $H^1$ yields
  \begin{align*}
    \|\theta^\kappa_t\|_{L^2} &\leq C\|\theta^\kappa_t\|_{H^{-\delta}}^{1/(1+\delta)} \|\theta^\kappa_t\|_{H^1}^{\delta/(1+\delta)}\\
                  &\leq C\|\theta_0\|_{L^2}\kappa^{-1/2}\widetilde D_0^{1/(1+\delta)}D_1t^2
                  \exp\left(\left(\frac{C}{1+\delta}+\frac{1}{2}\right)V(u_0)\right)e^{-\gamma t/(1+\delta)}\\
                  &\leq \|\theta_0\|_{L^2} \widehat D \kappa^{-1/2} \exp\left(\left(C + \frac{1}{2}\right)V(u_0)\right) e^{-\gamma t}.
  \end{align*}
  where we choose $\gamma$ even smaller to cancel out the factor $t^2$, and define $\widehat D \coloneq C\widetilde D_0^{1/(1+\delta)}D_1$.
  We conclude that $\mathbf E[\widehat D]\leq C$ by H\"older's inequality, since $\mathbf E[D_1^{1 + 1/\delta}] < \infty$.
  Finally, $L^2$ contraction gives, for $t\geq 2$,
  \begin{align*}
    \|\theta^\kappa_t\|_{L^2}
    &\leq \|\theta_0\|_{L^2}\min\left\{1,\widehat D\kappa^{-1/2}\exp\left(\left(C+\frac12\right)V(u_0)\right)e^{-\gamma t}\right\}\\
    &\leq \|\theta_0\|_{L^2}\widehat D^\delta\kappa^{-\delta/2}
    \exp\left(\delta\left(C+\frac12\right)V(u_0)\right)e^{-\delta\gamma t}.
  \end{align*}
  Taking $D=e^{2\delta\gamma}+\widehat D^\delta$ and adjusting $\gamma,\chi$ yields the conclusion for all $t\geq 0$, with $\mathbf E[D]\leq C$.
\end{proof}

\section{Thermally driven weakly active scalars}
\label{s:thermal_flows}
In this section, we are interested in the following hyperviscous Boussinesq system on the two-dimensional torus $\T^2$ with small coupling:
\begin{equation}\label{eq:eps-hNSE}
  \begin{cases}
    u^\varepsilon = \nabla^\perp \Delta^{-1} \omega^\varepsilon\\
    \partial_t \omega^\varepsilon + u^\varepsilon \cdot \nabla \omega^\varepsilon +\nu \Delta^{2} \omega^\varepsilon = \dot{F} + \varepsilon \Gamma \theta^\varepsilon &\quad \text{for } t > 0,\\
    \partial_t \theta^\varepsilon + u^\varepsilon \cdot \nabla \theta^\varepsilon = 0 &\quad \text{for } t > 0,\\
    \omega^\varepsilon = \omega_0, \quad \theta^\varepsilon = \theta_0 &\quad \text{for } t = 0,
  \end{cases}
\end{equation}
where $\Gamma: L^2_0(\T^2)\to H^{-1}_0(\T^2)$ is a bounded linear operator satisfying also $\Gamma: H^{-1}(\T^2)\to H^{-2}(\T^2)$, $\nu>0$, $\varepsilon \in \R$ is the buoyancy coefficient, the forcing term $F$ is given by
\begin{equation}
  \label{eq:F-def}
  F(t, x) \coloneq \sum_{k \in \Z^2 \setminus \{0\}} \zeta_k e_k(x) W^k_t,
\end{equation}
and $\zeta_k \in [c|k|^{-\ell}, C|k|^{-\ell}]$ for some fixed $\ell > 4$, $\{W^k\}$ is a family of independent standard Brownian motions, $\omega_0 \in L^2_0(\T^2)$, and
\begin{equation}
  \label{eq:ek-def}
  e_k(x) \coloneq
  \begin{cases}
    \sin(2 \pi k \cdot x) &\quad \text{if $k > (0, 0)$ lexicographically,}\\
    \cos(2 \pi k \cdot x) &\quad \text{otherwise.}
  \end{cases}
\end{equation}
The main result of this section reads as follows. In this section we allow constants $C, c > 0$ to depend on $\nu$, $\Gamma$, and $\{\zeta_k\}$ but on nothing else unless otherwise specified.

\begin{theorem}[Negative regularity mixing for a weakly active scalar]
\label{t:hNSE-mixing}
Under the above assumptions, there are $\gamma, \sigma, \varepsilon_0, \delta > 0$ such that for all $\varepsilon \in [-\varepsilon_0, \varepsilon_0]$ and $t \geq 0$,  \begin{equation}
    \label{eq:hNSE-mix}
    \E \|\theta^\varepsilon_t\|_{H^{-\delta}(\T^2)} \leq Ce^{- \gamma t} e^{\sigma\|\omega_0\|_{L^2}^2 }  \|\theta_0\|_{H^{-\delta}(\T^2)}
  \end{equation}
for all $\omega_0 \in L^2_0(\T^2)$ and $\theta_0\in L_0^2(\T^2)$ such that $\|\theta_0\|_{L^\infty(\T^2)} \leq 1$.
\end{theorem}

The proof of the above is postponed to Subsection~\ref{ss:proof_hNSE-mixing}.

\begin{remark}[Further applications to thermally driven flows]
The arguments used to prove~\eqref{eq:hNSE-mix} via the abstract results of Theorems~\ref{thm:mainStability} and~\ref{thm:shortTimeToLong} apply to a wider class of problems. In particular, we believe that the hyperviscosity in~\eqref{eq:eps-hNSE} is not essential, and more refined $L^p$-estimates with $p\in (1,\infty)$ can be used at least in the physically relevant case $\Gamma=-\partial_1$. For simplicity of exposition, we do not pursue this here.
For instance, the negative-mixing estimate of~\eqref{eq:hNSE-mix} also holds for solutions to the Stokes-Boussinesq system on the $d$-dimensional torus in all dimensions $d\geq 2$:
    \begin{equation}
\label{eq:stokes_boussinesq}
\left\{
\begin{aligned}
&\partial_t u^\varepsilon =-\nabla p^\varepsilon +\Delta u^\varepsilon-\varepsilon \theta^\varepsilon e_d+ \dot{F} &\quad \text{for } t > 0,\\
&\nabla \cdot u^\varepsilon=0&\quad \text{for } t > 0,\\
&\partial_t \theta^\varepsilon+ u^\varepsilon\cdot \nabla \theta^\varepsilon=0&\quad \text{for } t > 0,\\
&u^\varepsilon=u_0, \qquad
\theta^\varepsilon=\theta_0&\quad \text{for } t =0.
\end{aligned}
\right.
\end{equation}
where $e_d=(0,\dots,0,1)$, and $F$ is the standard divergence-free $d$-dimensional analogue of~\eqref{eq:F-def}, with Fourier coefficients $\eqsim |k|^{-\ell}$ for some $\ell\gg 1$. We leave the details to the interested reader.
\end{remark}

Before applying Theorems~\ref{thm:mainStability} and~\ref{thm:shortTimeToLong}, we describe how the hyperviscous Boussinesq system~\eqref{eq:eps-hNSE} fits into the abstract framework of Section~\ref{sec:abstract_results}.
Fix $T>1$, to be chosen sufficiently large below, and let
\begin{equation}
\label{eq:def_U_Hyper}
U=L^2_0(\T^2)\times \mathcal{B} \qquad \text{ where }\qquad
\mathcal{B}\coloneq \{\theta\in L^2_0(\T^2)\,:\, \|\theta\|_{L^\infty(\T^2)}\leq 1\},
\end{equation}
where $\mathcal{B}$ is endowed with the $H^{-1}(\T^2)$-topology. Moreover, we define
\begin{equation}
\label{eq:def_Z_Hyper}
Z\coloneq \{ z=(\omega,\theta) \in C([0,T];U)\,:\,  \om\in L^2(0,T;H^2)\},
\end{equation}
equipped with the natural intersection topology of $C([0,T];U)\cap L^2(0,T;H^2)$ on the vorticity component. We anticipate that the second condition on $\omega$ is used below to ensure the unique solvability of the associated ODE flow.

For $a=(\omega_0,\theta_0)\in U$ and $\varepsilon\in[-1,1]$, denote by $z_a^\varepsilon$ the trajectory of~\eqref{eq:eps-hNSE} on $[0,T]$ starting from $a$.
The maximum principle implies that, almost surely,
\begin{equation}
\label{eq:uniform_bound_theta_thermal_fluids}
\|\theta_t^\varepsilon\|_{L^\infty}\leq 1 \ \ \text{ for all }t\geq 0.
\end{equation}
Moreover, standard energy estimates established below yield
\[
\omega^\varepsilon \in C([0,T];L^2_0(\T^2)) \cap L^2(0,T;H^2(\T^2))
\]
a.s.\ (see Subsection~\ref{ss:apriori_estimates_Boussinesq}). Together with a Galerkin approximation and the initial-data stability estimates of Lemma~\ref{lem:hNSE-initial-data-stability} below, these bounds give global well-posedness and measurable dependence on the initial data. In particular, $z_a^\varepsilon\in Z$ almost surely.

For $z\in Z$, let $a_z=z(T)\in U$. We define the transition kernels $P_1,Q_1:Z\to\Prob(Z)$ by
$$
P_1(z,A)=\P(z^0_{a_z}\in A),
\qquad
Q_1(z,A)=\P(z^\varepsilon_{a_z}\in A),
$$
for every Borel set $A\subseteq Z$. The above describes successive trajectories on subsequent time intervals of length $T$. One can readily check that they form a Markov chain on $Z$, whose transition kernel is $P_1$ or $Q_1$, respectively. Next, for $z=(\omega,\theta)\in Z$, set
$
u_z=\nabla^\perp\Delta^{-1}\omega
$
and denote by $\Phi_t^z$ the flow solving
$$
\tfrac{\dd}{\dd t}\Phi_t^z(x)
=u_z(t,\Phi_t^z(x)),
\qquad
\Phi_0^z(x)=x.
$$
Since
$
u_z\in L^2(0,T;H^3(\T^2))
\subseteq L^1(0,T;C^{1,\eta}(\T^2))
$
for every $\eta\in(0,1)$, this ODE generates a unique volume-preserving $C^1$-diffeomorphism. We may therefore define the measurable maps
$$
\varphi_z=\psi_z:=\Phi_T^z \in\Diff^1(\T^2).
$$
For fixed $a=(\omega_0,\theta_0)\in U$ and $\varepsilon\in \R$, we use the shorthand
\begin{equation}\label{eq:phi_t_eps_defn}
\varphi_t^\varepsilon\coloneq\Phi_t^{z_a^\varepsilon}, \text{ for }t\geq 0.
\end{equation}
Finally, for the Wasserstein smallness condition~\eqref{eq:smallness_condition_Wasmodified}, we employ  Remark~\ref{r:Was_condition_checked}.
By choosing the coupling $
\Pi_z =\operatorname{Law} (z^0_{a_z},z^\varepsilon_{a_z}) \in\Gamma(P_1(z),Q_1(z)),
$
in~\eqref{eq:condition_smallness_reduced},  from the above choice, we obtain
\begin{equation}
\label{eq:Was_difference_epsilon_thermal_fluids}
\Wasl_{d\wedge 1,m}[\mu_z,\nu_z] \leq \E\Big[ (1+\BiLip(\varphi_{z^0_{a_z}}) +\BiLip(\varphi_{z^\varepsilon_{a_z}}))^m
(\dist^1(\varphi_{z^0_{a_z}},
\varphi_{z^\varepsilon_{a_z}})\wedge 1)
\Big],
\end{equation}
where $m=p(1+\delta)+1$ and $p\in (0,2)$ and $\delta>0$ are specified later. This quantity will be estimated in Subsection~\ref{ss:stability_in_epsilon_thermal}.

\subsection{A priori estimates}
\label{ss:apriori_estimates_Boussinesq}
In this subsection, we prove uniform in $\varepsilon\in [-1,1]$ estimates for the following (logarithm of the) Lyapunov function for~\eqref{eq:eps-hNSE}:
\begin{equation}
  \label{eq:V-def-hNSE}
  V_{\sigma}(\omega) \coloneq \sigma \|\omega\|_{L^2}^2,
\end{equation}
where $\sigma > 0$. Here, we follow the argument of~\cite[Appendix~A.1]{CR26} (see also~\cite{BBPS22}). Below we use the uniform bound~\eqref{eq:uniform_bound_theta_thermal_fluids} for the temperature field.

\begin{proposition}\label{prop:V-hNSE}
For each $p, T \geq 1$ and $\delta > 0$, there are $\alpha > 1$ and $r, \sigma, C > 0$ with
\begin{align}
  \label{eq:V-contraction-hNSE}
  \E \left[\exp\left(\alpha V_{\sigma}(\omega^\varepsilon_T)\right) \right] &\leq C \exp(V_{\sigma}(\omega_0)),\\
  \label{eq:V-bounded-hNSE}
  \E \left[\exp\left( \int_0^T r\|\nabla u^\varepsilon_s\|_{C^r}^{\alpha} \dif s\right) \right] &\leq C \exp(\delta V_{\sigma}(\omega_0)),\\
  \label{eq:phi-bounded-hNSE}
  \E \left[\sup_{t \in [0, T]} \BiLip(\varphi^\varepsilon_t)^p\right] &\leq C \exp(\delta V_{\sigma}(\omega_0)),
\end{align}
for all $\varepsilon \in [-1, 1]$.
\end{proposition}

\begin{proof}
  Lemma~\ref{lem:L2-contraction-hNSE} immediately gives~\eqref{eq:V-contraction-hNSE}.
  Combine Lemma~\ref{lem:L2-contraction-hNSE}, Gagliardo--Nirenberg and Morrey's inequality for~\eqref{eq:V-bounded-hNSE}.
Finally,~\eqref{eq:phi-bounded-hNSE} follows from Lemma~\ref{lem:velocityToFlowMap} and~\eqref{eq:V-bounded-hNSE}.
\end{proof}

\begin{lemma}\label{lem:L2-contraction-hNSE}
For each $T > 0$, there is $\eta_0 > 0$ such that for all $\eta \in [0, \eta_0]$,
\begin{equation}\label{eq:L2-contraction-hNSE}
  \E \left[ \sup_{t \in [0, T]}\exp\left(\eta \left(e^{\nu t} \|\omega^\varepsilon_t\|_{L^2}^2 + c \int_0^t \|\Delta \omega^\varepsilon_s\|_{L^2}^2 \dif s\right)\right)\right] \leq C\exp(\eta \|\omega_0\|_{L^2}^2)
\end{equation}
for all $\varepsilon \in [-1, 1]$.
\end{lemma}

\begin{proof}
  Applying It\^o's lemma to $e^{\nu t}\|\omega^\varepsilon_t\|_{L^2}^2$,
  \begin{align*}
    \dif \left(e^{\nu t} \|\omega^\varepsilon_t\|_{L^2}^2\right) &= e^{\nu t} \left(\left(\nu \|\omega^\varepsilon_t\|_{L^2}^2 + \sum_{k \in \Z^2 \setminus \{0\}} |\zeta_k|^2\|e_k\|_{L^2}^2\right) \dif t + 2 \langle \dif \omega^\varepsilon_t, \omega^\varepsilon_t \rangle\right)\\
&= e^{\nu t}\Bigg( \left(\nu \|\omega^\varepsilon_t\|_{L^2}^2 + \sum_k |\zeta_k|^2\|e_k\|_{L^2}^2 + 2 \left( \langle \omega^\varepsilon_t,- u^\varepsilon_t \cdot \nabla \omega^\varepsilon_t + \varepsilon \Gamma \theta^\varepsilon_t \rangle - \nu \|\Delta \omega^\varepsilon_t\|_{L^2}^2\right)\right)\dif t \\
&\qquad + 2\langle \omega^\varepsilon_t, \dif F_t \rangle\Bigg)\\
&\leq e^{\nu t} \left(\left(C - \frac{\nu}{2}\|\Delta \omega^\varepsilon_t\|_{L^2}^2\right)\dif t + 2\langle \omega^\varepsilon_t, \dif F_t \rangle\right),
  \end{align*}
  where we use the fact that $\|\Delta^{-1} \Gamma \theta_t^\varepsilon\|_{L^2} \leq C\|\theta_t^\varepsilon\|_{H^{-1}}\leq C$.
  Integrating from $0$ to $t$, we have
  \begin{align*}
    e^{\nu t} \|\omega^\varepsilon_t\|_{L^2}^2 - \|\omega_0\|_{L^2}^2 + \int_0^t \frac{\nu}{4}\|\Delta \omega^\varepsilon_s\|_{L^2}^2 \dif s &\leq \int_0^t e^{\nu s}\left(C - \frac{\nu}{4} \|\omega^\varepsilon_s\|_{L^2}^2\right) \dif s + 2\sum_k \int_0^t e^{\nu s} \zeta_k \langle \omega^\varepsilon_s, e_k\rangle \dif W^k_s\\
                                                                                        &\leq Ce^{\nu t} - ce^{-\nu t}[M]_t + M_t,
  \end{align*}
  where
  \begin{equation*}
    M_t \coloneq 2\sum_k \int_0^t e^{\nu s} \zeta_k \langle \omega^\varepsilon_s, e_k\rangle \dif W^k_s
  \end{equation*}
  is a martingale.
  Indeed, since $|\zeta_k| \approx |k|^{-\ell}$ and $\ell > 1$, the quadratic variation of $M$ satisfies $\frac{\dif}{\dif s}[M]_s \leq Ce^{2\nu s}\|\omega^\varepsilon_s\|_{L^2}^2$. For $\eta \in (0, ce^{-\nu T})$, we deduce
  \begin{equation*}
    \eta \left(e^{\nu t}\|\omega^\varepsilon_t\|_{L^2}^2 - \|\omega_0\|_{L^2}^2 + \int_0^t \frac{\nu}{4} \|\Delta \omega^\varepsilon_s\|_{L^2}^2 \dif s\right) \leq C - [\eta M]_t + \eta M_t.
  \end{equation*}
and we conclude~\eqref{eq:L2-contraction-hNSE} by applying Doob's maximal inequality to the nonnegative supermartingale $\exp(2\eta M_t-2[\eta M]_t)$ and then taking square roots.
\end{proof}

\begin{remark}
\label{rem:large-alpha-hNSE}
From the above proofs, it follows that, for any $A>1$, choosing $T$ sufficiently large,
the estimate~\eqref{eq:V-contraction-hNSE} holds with $A$ in place
of $\alpha$. Indeed, fix $A>1$
and choose $T$ sufficiently large so that
$Ae^{-\nu T}\leq 1$. Then, by
Lemma~\ref{lem:L2-contraction-hNSE}, choosing
$\sigma>0$ sufficiently small, we have
$$
\E [\exp(A V_\sigma(\omega^\varepsilon_T))]
\leq C \exp( Ae^{-\nu T}V_\sigma(\omega_0))
\leq C \exp(V_\sigma(\omega_0)).
$$
\end{remark}

By Galerkin approximation and Proposition~\ref{prop:V-hNSE} we see that solutions to~\eqref{eq:eps-hNSE} exist for all time almost surely.
We now turn to stability of~\eqref{eq:eps-hNSE} in the initial data.

\begin{lemma}
\label{lem:hNSE-initial-data-stability}
  Let $(u^\varepsilon, \omega^\varepsilon, \theta^\varepsilon)$ and $(\tilde{u}^\varepsilon, \tilde{\omega}^\varepsilon, \tilde{\theta}^\varepsilon)$ be solutions to~\eqref{eq:eps-hNSE} with initial data $(u_0, \omega_0, \theta_0)$ and $(\tilde{u}_0, \tilde{\omega}_0, \tilde{\theta}_0)$ respectively and identical forcing, with $(\omega_0,\theta_0),(\tilde{\omega}_0,\tilde{\theta}_0)\in U$ (see~\eqref{eq:def_U_Hyper}).
  Then for every $T, K > 0$ there are $C > 0$ and $\sigma > 0$ such that
  \begin{equation}
    \label{eq:eps-hNSE-initial-data-stability}
    \E \sup_{t \in [0, T]} \left(\log \left(1 + \frac{\|\omega^\varepsilon_t - \tilde{\omega}_t^\varepsilon\|_{L^2}^2 + \|\theta^\varepsilon_t - \tilde{\theta}^\varepsilon_t\|_{H^{-1}}^2}{\|\omega_0 - \tilde{\omega}_0\|_{L^2}^2 + \|\theta_0 - \tilde{\theta}_0\|_{H^{-1}}^2}\right)\right)^K \leq C\exp\left(V_{\sigma}(\omega_0)\right),
  \end{equation}
  where we take the convention that $\frac{0}{0} = 1$ and $\frac{a}{0} = +\infty$ for any $a > 0$.
\end{lemma}

\begin{proof}
  Write $\xi \coloneq \omega^\varepsilon - \tilde{\omega}^\varepsilon$ and $\chi \coloneq \theta^\varepsilon - \tilde{\theta}^\varepsilon$.
  Then $\xi$ solves
  \begin{equation}
    \label{eq:difference-hNSE-stability-initial}
    \partial_t \xi + u^\varepsilon \cdot \nabla \omega^\varepsilon - \tilde{u}^\varepsilon\cdot \nabla \tilde{\omega}^\varepsilon + \nu (-\Delta)^2 \xi =  \varepsilon \Gamma \chi,\\
  \end{equation}
  so integrating against $\xi$ yields
  \begin{align*}
    \frac{\dif}{\dif s} \|\xi_s\|_{L^2}^2 + 2\nu \|\Delta \xi_s\|_{L^2}^2 &= 2\int \xi_s \left( \tilde{u}^\varepsilon_s\cdot \nabla \tilde{\omega}^\varepsilon_s - u^\varepsilon_s \cdot \nabla \omega^\varepsilon_s +  \varepsilon \Gamma \chi_s \right)\\
                                                           &= 2\int \xi_s \left( (\tilde{u}^\varepsilon_s - u^\varepsilon_s)\cdot \nabla \omega^\varepsilon_s +  \varepsilon \Gamma \chi_s \right)\\
                                                           &\leq 2\|\xi_s\|_{L^2}\|\tilde{u}^\varepsilon_s - u^\varepsilon_s\|_{L^4}\|\nabla \omega^\varepsilon_s\|_{L^4}  + C|\varepsilon| \|\Delta \xi_s\|_{L^2}\|\chi_s\|_{H^{-1}}\\
                                                           &\leq C\left(\|\xi_s\|_{L^2}^2 \|\nabla \omega^\varepsilon_s\|_{L^4} + \varepsilon^2 \|\chi_s\|_{H^{-1}}^2\right) + \nu \|\Delta \xi_s\|_{L^2}^2\\
                                                           &\leq C\left(\|\xi_s\|_{L^2}^2 \|\Delta \omega^\varepsilon_s\|_{L^2} + \varepsilon^2 \|\chi_s\|_{H^{-1}}^2\right) + \nu \|\Delta \xi_s\|_{L^2}^2.
  \end{align*}
  By Gr\"onwall,
  \begin{equation}\label{eq:hxi-T-Gronwall}
    \|\xi_t\|_{L^2}^2 \leq \left(\|\xi_0\|_{L^2}^2 + C\varepsilon^2\int_0^t \|\chi_s\|_{H^{-1}}^2 \dif s\right)\exp\left(C\int_0^t \|\Delta \omega^\varepsilon_s\|_{L^2} \dif s \right).
  \end{equation}

  On the other hand, the $H^1$ transport estimate and duality give
  \begin{align*}
    \frac{\dif}{\dif s} \|\chi_s\|_{H^{-1}}^2 &\leq C \left(\|\chi_s\|_{H^{-1}}^2\|\nabla u^\varepsilon_s\|_{L^\infty} + \|\chi_s\|_{H^{-1}}\|(u^\varepsilon_s - \tilde{u}^\varepsilon_s)  \tilde{\theta}^\varepsilon_s\|_{L^2}\right)\\
                                         &\leq C\|\chi_s\|_{H^{-1}}^2(1 + \|\nabla u^\varepsilon_s\|_{L^\infty}) + C\|(u^\varepsilon_s - \tilde{u}^\varepsilon_s)  \tilde{\theta}^\varepsilon_s\|_{L^2}^2\\
                                         &\leq C\|\chi_s\|_{H^{-1}}^2(1 + \|\nabla u^\varepsilon_s\|_{L^\infty}) + C\|\xi_s\|_{L^2}^2 \|\tilde{\theta}^\varepsilon_s\|_{L^4}^2.
  \end{align*}
  and therefore
  \begin{equation}
    \label{eq:hchi-T-Gronwall}
    \|\chi_t\|_{H^{-1}}^2 \leq \left(\|\chi_0\|_{H^{-1}}^2 + C\int_0^t \|\xi_s\|_{L^2}^2 \dif s\right)\exp\left(C\int_0^t (1 + \|\nabla u^\varepsilon_s\|_{L^\infty}) \dif s\right).
  \end{equation}
  Combining~\eqref{eq:hxi-T-Gronwall} and~\eqref{eq:hchi-T-Gronwall}, we have
  \begin{multline*}
    \|\xi_t\|_{L^2}^2 + \|\chi_t\|_{H^{-1}}^2\leq \left(\|\xi_0\|_{L^2}^2 + \|\chi_0\|_{H^{-1}}^2 + C\int_0^t \left(\|\chi_s\|_{H^{-1}}^2 + \|\xi_s\|_{L^2}^2\right) \dif s\right)\\ {}\times\exp\left(C\int_0^t \left(1 + \|\Delta \omega^\varepsilon_s\|_{L^2}\right) \dif s \right)
  \end{multline*}
  and so by one more application of Gr\"onwall,
  \begin{align*}
    \sup_{t \in [0, T]}\big( \|\xi_t\|_{L^2}^2 + \|\chi_t\|_{H^{-1}}^2\big) \leq \left(\|\xi_0\|_{L^2}^2 + \|\chi_0\|_{H^{-1}}^2\right)&\exp\left(C\int_0^T \left(1 + \|\Delta \omega^\varepsilon_s\|_{L^2}\right) \dif s \right)\\
    \cdot\exp\Bigg(&\exp\left(C\int_0^T \left(1 + \|\Delta \omega^\varepsilon_s\|_{L^2}\right) \dif s \right)\Bigg).
  \end{align*}
We conclude with Lemma~\ref{lem:L2-contraction-hNSE}.
\end{proof}

\subsection{Stability in \texorpdfstring{$\varepsilon$}{epsilon}}
\label{ss:stability_in_epsilon_thermal}
In this subsection, we prove that the diffeomorphisms $\varphi_{z^0_{a_z}}$ and $\varphi_{z^\varepsilon_{a_z}}$ appearing in~\eqref{eq:Was_difference_epsilon_thermal_fluids} are close quantitatively in $\varepsilon$. We begin with suitable stability estimates for the corresponding vector fields.

\begin{lemma}\label{lem:homega-eps-stability}
  Let $(u^\varepsilon, \omega^\varepsilon, \theta^\varepsilon)$ and $(u^0, \omega^0, \theta^0)$ be solutions to~\eqref{eq:eps-hNSE} with identical forcing and initial data $(u_0, \omega_0, \theta_0)$ with $(\omega_0,\theta_0)\in U$ (see~\eqref{eq:def_U_Hyper}).  Then for each $T, K > 0$ there are $C, \sigma > 0$ such that
  \begin{equation}
    \label{eq:eps-hNSE-stability}
    \E \sup_{t \in [0, T]} \|\omega^\varepsilon_t - \omega^0_t\|_{L^2}^K \leq C|\varepsilon|^K\exp\left(V_{\sigma}(\omega_0)\right).
  \end{equation}
\end{lemma}
\begin{proof}
  Write $\xi \coloneq \omega^\varepsilon - \omega^0$ and observe that
  \begin{equation*}
    \partial_t \xi + u^\varepsilon \cdot \nabla \omega^\varepsilon - u^0 \cdot \nabla \omega^0 + \nu (-\Delta)^2 \xi = \varepsilon \Gamma \theta^\varepsilon.
  \end{equation*}
  Integrating against $2 \xi$, we have
  \begin{align*}
    \frac{\dif}{\dif s} \|\xi_s\|_{L^2}^2 + 2\nu \|\Delta \xi_s\|_{L^2}^2 &\leq 2\int \xi_s (u^0_s \cdot \nabla \omega^0_s - u^\varepsilon_s \cdot \nabla \omega^\varepsilon_s ) + 2\varepsilon\langle\Gamma\theta^\varepsilon_s,\xi_s\rangle\\
                                                           &\leq 2 \int \xi_s \left((u^0_s - u^\varepsilon_s) \cdot \nabla \omega^0_s \right) + 2\varepsilon\langle\Gamma\theta^\varepsilon_s,\xi_s\rangle\\
                                                           &\leq 2 \|\xi_s\|_{L^2} \|u^0_s - u^\varepsilon_s\|_{L^4} \|\nabla \omega^0_s\|_{L^4} +C |\varepsilon| \|\Delta \xi_s\|_{L^2}\|\theta^\varepsilon_s\|_{H^{-1}}\\
                                                           &\leq 2 \|\xi_s\|_{L^2} \|\xi_s\|_{L^2} \|\nabla \omega^0_s\|_{L^4} + C\varepsilon^2 \|\theta^\varepsilon_s\|_{H^{-1}}^2 + \nu\|\Delta \xi_s\|_{L^2}^2\\
                                                           &\leq C\|\xi_s\|_{L^2}^2 \|\Delta \omega^0_s\|_{L^2} + C\varepsilon^2 \|\theta^\varepsilon_s\|_{H^{-1}}^2 + \nu\|\Delta \xi_s\|_{L^2}^2.
  \end{align*}
  and by Gr\"onwall
  \begin{equation*}
    \sup_{t \in [0, T]} \|\xi_t\|_{L^2}^2 \leq C\varepsilon^2\exp\left(C\int_0^T \|\Delta \omega^0_t\|_{L^2} \dif t\right).
  \end{equation*}
  We conclude with Lemma~\ref{lem:L2-contraction-hNSE}.
\end{proof}

\begin{lemma}\label{lem:hflow-eps-stability}
For each $p,a,T,\sigma>0$ there are $C,c>0$ such that, for all $\varepsilon\in[-1,1]$,
\begin{equation*}
\E \sup_{t \in [0, T]} \dist^1(\varphi^\varepsilon_t, \varphi^0_t)^p \leq C|\varepsilon|^c \exp(aV_{\sigma}(\omega_0)).  \end{equation*}
\end{lemma}
\begin{proof}
  Using Lemma~\ref{lem:velocityToFlowMap}, Proposition~\ref{prop:V-hNSE}, Lemma~\ref{lem:L2-contraction-hNSE} and Lemma~\ref{lem:homega-eps-stability}, for fixed small $\delta > 0$ we have
  \begin{align*}
    \E \sup_{t \in [0, T]} \dist^1(\varphi^\varepsilon_t, \varphi^0_t)^p &\leq C \E \left[\left(\int_0^T\|u^\varepsilon_t - u^0_t\|_{L^\infty} \dif t\right)^{\frac{p\delta}{1+\delta}} \exp\left(C \int_0^T \left(\|u^\varepsilon_t\|_{C^{1+\delta}} + \|u^0_t\|_{C^{1+\delta}}\right) \dif t\right)\right]\\
                           &\leq C \exp(V_{\sigma}(\omega_0))\E \left[\left(\int_0^T\|u^\varepsilon_t - u^0_t\|_{L^\infty} \dif t\right)^p\right]^{\frac{\delta}{1+\delta}}\\
                           &\leq C \exp(V_{\sigma}(\omega_0))\E \left[\left(\int_0^T\|\omega^\varepsilon_t - \omega^0_t\|_{L^2}^{\delta}\|\omega^\varepsilon_t - \omega^0_t\|_{H^1}^{1-\delta} \dif t\right)^p\right]^{\frac{\delta}{1+\delta}}\\
                           &\leq C |\varepsilon|^{\frac{\delta^2 p}{1+\delta}}\exp(2V_{\sigma}(\omega_0)) \E \left[\left(\int_0^T \left(\|\omega^\varepsilon_t\|_{H^1} + \|\omega^0_t\|_{H^1}\right)^{(1-\delta)} \dif t\right)^{2p}\right]^{\frac{\delta}{2(1+\delta)}}\\
                           &\leq C |\varepsilon|^{\frac{\delta^2 p}{1+\delta}}\exp(3V_{\sigma}(\omega_0)),
  \end{align*}
  where we fix $\delta \in (0,1)$, so that $H^{2-\delta}(\T^2)\hookrightarrow L^\infty(\T^2)$ and $H^3(\T^2)\hookrightarrow C^{1+\delta}(\T^2)$. The factor of $3$ in the exponential can be replaced by $a$ by redefining $\sigma$, since $3V_\sigma=aV_{3\sigma/a}$.
\end{proof}

\subsection{Proof of Theorem~\ref{t:hNSE-mixing}}
\label{ss:proof_hNSE-mixing}

\begin{proof}[Proof of Theorem~\ref{t:hNSE-mixing}]
Recall that $Z$ and $V_\sigma$ are defined in~\eqref{eq:def_Z_Hyper} and~\eqref{eq:V-def-hNSE}, respectively. Next, extend $V_\sigma$ to $Z$ by setting $V(z)\coloneq V_\sigma(\omega(T))$ for $z=(\omega,\theta)\in Z$. Fix $p\in(1,2)$ and $q=(1+p)/2$.
By Proposition~\ref{prop:V-hNSE},
Remark~\ref{rem:large-alpha-hNSE}
and~\cite[Theorem 1.8]{BFPS24},
choosing $T$ sufficiently large and $\sigma>0$ suitably,
Assumptions~\ref{ass:negative_regularity_mixing}
and~\ref{ass:super_lyapunov} hold for $P_1$,
with constants $\alpha,K,\beta$ satisfying
$$
\frac{pq}{p-q}\leq\alpha-1 \quad \text{ and } \beta K^{1/(\alpha-1)}<\frac{1}{2}.
$$
Strictly speaking, the (logarithm of the) super--Lyapunov function in~\cite[Theorem~1.8 and Section~2.2]{BFPS24} is of the form
\[
W_\sigma(\omega)\coloneq M \log (1+\|\omega\|_{H^s}) + V_\sigma(\omega),
\]
where $M, s>0$ are fixed and $\sigma>0$ is sufficiently small.
However, a standard argument allows us to replace this weight with $V_\sigma$.
Indeed, we simply wait for time $1$ (for $\omega$ to smooth) before applying the contraction estimate.
One readily checks that, for $\sigma'>0$ sufficiently small depending on $\sigma$ and $\alpha$,
\begin{align*}
\E \exp(\alpha W_{\sigma'}(\omega^0_1)) &\leq C e^{V_\sigma(\omega_0)},\\
\E \BiLip(\varphi^0_1)^r &\leq C_r e^{V_\sigma(\omega_0)},
\end{align*}
for any $r>0$. Since $\alpha>p$, H\"older's inequality with exponents $\alpha/p$ and $\alpha/(\alpha-p)$ gives
\[
\E\left[e^{pW_{\sigma'}(\omega^0_1)}\BiLip(\varphi^0_1)^{p(1+\delta)}\right]\leq C e^{V_\sigma(\omega_0)},
\]
and conditioning the contraction estimate at time $1$ yields the required contraction with weight $V_\sigma$ for sufficiently large $T$.

Moreover, by~\eqref{eq:Was_difference_epsilon_thermal_fluids}, Proposition~\ref{prop:V-hNSE} and Lemma~\ref{lem:hflow-eps-stability}, there are $C,c>0$ such that
$$
\Wasl_{d\wedge 1,m}[\mu_z,\nu_z]\leq C|\varepsilon|^c e^{pV(z)}.
$$
From Theorem~\ref{thm:mainStability} with $\psi=\varphi$, it follows that Assumption~\ref{ass:negative_regularity_mixing} holds for $Q_1$ with $\beta'K^{1/(\alpha-1)}<1$, provided $|\varepsilon|\leq\varepsilon_0$ is sufficiently small. Together with~\eqref{eq:V-contraction-hNSE}, Theorem~\ref{thm:shortTimeToLong} therefore gives
$$
(\E\|\theta^\varepsilon_{nT}\|_{H^{-\delta}}^q)^{1/q}\leq C e^{-\gamma n} e^{\sigma\|\omega_0\|_{L^2}^2} \|\theta_0\|_{H^{-\delta}}  \ \ \text{ for all } n\in\N,
$$
after increasing $\sigma>0$ on the right-hand side if necessary.
Finally, we extend the above estimate to all times $t\geq0$. To this end, write $n=\floor{t/T}$ and $q'=q/(q-1)$ (recall that $q=(1+p)/2>1$).
Let $\varphi^\varepsilon_{s,t}$ be the Lagrangian flow map from time $s$ to time $t$ induced by the process $(\omega^\varepsilon_t)_{t\geq 0}$.
By~\eqref{eq:H-delta-interp} and H\"older's inequality,
\begin{align*}
\E\|\theta^\varepsilon_t\|_{H^{-\delta}}&=\E\|(\varphi^\varepsilon_{nT,t})_*\theta^\varepsilon_{nT}\|_{H^{-\delta}}\\
&\leq \E\Big[ \|\theta^\varepsilon_{nT}\|_{H^{-\delta}}\BiLip(\varphi^\varepsilon_{nT,t})^{1+\delta}\Big]\\
&\leq\big(\E\|\theta^\varepsilon_{nT}\|_{H^{-\delta}}^q\big)^{1/q}
\Big(\E\BiLip(\varphi^\varepsilon_{nT,t})^{q'(1+\delta)}\big)^{1/q'}\\
&\leq C e^{-\gamma n}  e^{\sigma\|\omega_0\|_{L^2}^2} \|\theta_0\|_{H^{-\delta}},
\end{align*}
where the last inequality follows from~\eqref{eq:V-contraction-hNSE},~\eqref{eq:phi-bounded-hNSE} and the Markov property.  The conclusion~\eqref{eq:hNSE-mix} follows by decreasing $\gamma$ and increasing $C$ and $\sigma$ if necessary.
\end{proof}

\section{Regularized forcing}\label{sec:regularized}

\subsection{Setup and assumptions}
In this section, we study the 2d Navier-Stokes equation with a continuous random body force on the domain $\mathbb T^2$. Written in vorticity form, this reads
\begin{equation}
  \label{eq:2dnse-smooth-noise}
  \begin{cases}
    \partial_t \omega^\varepsilon + u^\varepsilon \cdot \nabla \omega^\varepsilon - \nu \Delta \omega^\varepsilon = \xi^\varepsilon &\qquad \text{ for }x \in \mathbb T^2, t>0,\\
    \omega^\varepsilon= \omega_0, &\qquad\text{ for } x \in \mathbb T^2, t = 0,
  \end{cases}
\end{equation}
where $u^\varepsilon \coloneq \nabla^\perp \Delta^{-1} \omega^\varepsilon$.
We take the viscosity $\nu > 0$ to be a positive constant and consider deterministic initial data $\omega_0\in H^1_0$. We assume that the forcing is a mollification in time of the white-in-time, colored-in-space noise considered in~\cite{BBPS22}; that is, let $m,\ell > 0$ be constants satisfying
\[
\ell > 4, \qquad \max\{3, \ell -2\} < m < \ell -1.
\]
and let the random process $F_t \in H^m_0$ be given by
\begin{equation}\label{eq:reg-F-def}
F_t = \sum_{k \in  \mathbb Z^2 \setminus \{0\}} \beta_k e_k W^k_t, \qquad c |k |^{-\ell} \leq |\beta_k| \leq C|k|^{-\ell},
\end{equation}
where $\{W^k\}$ is a collection of independent standard (two-sided) Brownian motions
and $\{e_k\}_k$ is the real-valued Fourier basis from~\eqref{eq:ek-def}.
We normalize $F_0=0$.
Throughout this section, the initial noise history $(F_s)_{s\leq0}$ is sampled from its Brownian law, independently of the future increments.
We fix a nonnegative mollifier $\varrho \in C^1([0,\infty))$ satisfying
\begin{equation}\label{eq:reg-mollifier-assumptions}
  \int_0^\infty \varrho(s)\dif s=1,
  \qquad |\varrho(s)|+|\varrho'(s)|\leq C_\varrho e^{-as}
  \quad\text{for some }a,C_\varrho>0.
\end{equation}
Extend $\varrho$ by zero on $(-\infty,0)$ and define
\begin{equation}\label{eq:F-eps-def}
  F^\varepsilon_t \coloneq \int \varrho^\varepsilon(t-s) F_s \dif s, \qquad
  \xi^\varepsilon_t \coloneq \dot{F}^\varepsilon_t,
\end{equation}
where $\varrho^\varepsilon(x) = \varepsilon^{-1}\varrho(x/\varepsilon)$ for $\varepsilon>0$. The process $\xi^\varepsilon$ is continuous and adapted, hence predictable. At $\varepsilon=0$, we interpret $\varrho^\varepsilon = \delta_0$ and the forcing as $\dot F_t$.

\begin{remark}
  When $\varrho(t) = e^{-t}$, the regularized forcing is an Ornstein--Uhlenbeck process given by
  \begin{equation*}
    \eps \partial_t  \xi^\eps = - \xi^\eps +\xi^0.
  \end{equation*}
  More generally, mollifiers of the form $\varrho(t) = c t^ne^{-t}$ correspond to ``OU tower noise''; see~\cite[Remark 1.13]{BBPS22}.
\end{remark}

Our goal is to show that the negative-regularity mixing estimates also hold for small $\varepsilon > 0$.
It is not currently known whether there is a unique stationary measure for the $\omega^\varepsilon_t$-process when $\varepsilon > 0$, which is a key assumption involved in the proof of the negative-regularity mixing estimate.
Nevertheless, we show that negative-regularity mixing is stable under small perturbations of $\varepsilon$.

\begin{theorem}\label{thm:reg-noise-mixing}
Let $\omega^\eps$ solve~\eqref{eq:2dnse-smooth-noise} and let $\theta^\varepsilon_t$ solve
\begin{align*}
\partial_t\theta^\varepsilon_t + u^\eps_t \cdot \nabla_x \theta^\varepsilon_t = 0, \text{ for } x \in \mathbb T^2, t > 0,
\end{align*}
with $\theta^\varepsilon_0 = \theta_0 \in H^{-\delta}_0(\mathbb T^2)$.
There are $\eps_0,\delta,q,\gamma,C > 0$ with $\delta <1$ such that for all deterministic $(\omega_0,\theta_0) \in H^1_0 \times H^{-\delta}_0$,  $0 < \eps \leq \eps_0$, $t\geq 0 $, we have
\begin{align}\label{eq:reg_noise_neg_mixing}
\left(\mathbf E[\|\theta_t^\varepsilon\|_{H^{-\delta}}^q]\right)^\frac{1}{q} \leq C\exp\left(-\gamma t + C \|\omega_0\|_{H^1}^2\right) \|\theta_0\|_{H^{-\delta}}.
\end{align}
\end{theorem}

\subsection{Super-Lyapunov estimates}

For a continuous noise history with $F_0=0$, define
\begin{equation}\label{eq:reg-history-def}
  \|F\|_{\mathrm{hist}}\coloneq
  \sup_{s<0}\frac{e^{s/2}}{|s|^{1/3}}\|F_s\|_{H^m}.
\end{equation}
The exponent $1/3$ is chosen simply because it belongs to $(0, 1/2)$; it is useful that the history norm measures a bit of regularity of $F$ at time $0$.

We work on the measurable state space
\[
  Z\coloneq\big\{(\omega,F)\in C_{\mathrm{loc}}((-\infty,0];H^2_0\times H^m_0):F_0=0,\ \|F\|_{\mathrm{hist}}<\infty\big\},
\]
with the Borel $\sigma$-algebra inherited from the space of continuous histories.
Since $1/3<1/2$, Brownian histories belong to $Z$ in the noise coordinate almost surely, and Fernique's theorem gives $\E\exp(b\|F\|_{\mathrm{hist}}^2)<\infty$ for sufficiently small $b>0$. We abuse notation slightly and write $\|z\|_{\mathrm{hist}}$ to mean $\|F\|_{\mathrm{hist}}$ when $z = (\omega, F)$.

To describe the Markov transitions, condition on $z_0=\{\omega_s,F_s\}_{s\leq0}$, evolve $F$ by independent Brownian increments according to~\eqref{eq:reg-F-def}, and solve~\eqref{eq:2dnse-smooth-noise} from $\omega_0$. For $t\geq0$, recenter the noise history at its current value and set
\[
  z_t^\varepsilon\coloneq\{\omega^\varepsilon_{t+s},F_{t+s}-F_t\}_{s\leq0}.
\]
The histories form a Markov process whose time-$t$ transition kernel is denoted by $P_t$, for $t\geq0$. We write $\E_{z_0}$ for the corresponding conditional expectation. For $0\leq\varepsilon\leq\min\{1,a\}$, the convolution in~\eqref{eq:F-eps-def} is well-defined on $Z$ by~\eqref{eq:reg-mollifier-assumptions}.

Define
\begin{equation}\label{eq:reg-V-def}
  L_\rho(\omega)\coloneq\|\omega\|_{L^2}^2+\rho\|\nabla\omega\|_{L^2}^{2/3},
  \qquad V_{\sigma,\rho,b}(\omega,F)\coloneq\sigma L_\rho(\omega_0)+b\|F\|_{\mathrm{hist}}^2,
\end{equation}
where $\sigma,\rho,b>0$.

Our first goal is to prove that $V_{\sigma, \rho, b}$ is a super-Lyapunov function which controls enough regularity of the flow maps $\varphi_t^\varepsilon$, given by~\eqref{eq:phi_t_eps_defn}.

\begin{proposition}\label{prop:reg-V-NSE}
For each $p,T\geq1$ and $\delta>0$, there are $\varepsilon_0>0$, $\alpha>1$ and $r,\sigma,\rho,b,C>0$ such that
\begin{equation}\label{eq:reg-V-contraction-NSE}
  \E_{z_0}\exp\big(\alpha V_{\sigma,\rho,b}(z_T^\varepsilon)\big)
  \leq C\exp\big(V_{\sigma,\rho,b}(z_0)\big),
\end{equation}
\begin{equation}\label{eq:reg-V-bounded-NSE}
  \E_{z_0}\exp\left(\int_0^T r\|\nabla u_s^\varepsilon\|_{C^r}^{\alpha}\dif s\right)
  \leq C\exp\big(\delta V_{\sigma,\rho,b}(z_0)\big),
\end{equation}
and
\begin{equation}\label{eq:reg-phi-bounded-NSE}
  \E_{z_0}\left[\sup_{t\in[0,T]}\BiLip(\varphi_t^\varepsilon)^p\right]
  \leq C\exp\big(\delta V_{\sigma,\rho,b}(z_0)\big)
\end{equation}
for every $z_0\in Z$ and $\varepsilon\in[0,\varepsilon_0]$.
\end{proposition}

All constants below may depend on $T,\nu,m,(\beta_k)$ and $\varrho$, but not on $\varepsilon$ or $z_0$.

\begin{lemma}\label{lem:reg-noise-estimates}
There are $C_0,c,C>0$, depending only on $m,(\beta_k)$ and $\varrho$, such that
\begin{equation}\label{eq:reg-correction-tail}
  \mathbb P_{z_0}\left(
    \sup_{0\leq t\leq T}
    \left(\|F_t-F_t^\varepsilon\|_{H^m}
      -C_0\varepsilon^{1/3}e^{-at/\varepsilon}\|z_0\|_{\mathrm{hist}}\right)
    >\theta\right)
  \leq C(1+T/\varepsilon)e^{-c\theta^2/\varepsilon}
\end{equation}
for every $z_0\in Z$, $T,\theta>0$ and $0<\varepsilon\leq\min\{1,a\}$.
\end{lemma}
\begin{proof}
The conditional mean of $F_t-F_t^\varepsilon$ is
\[
  \E_{z_0}(F_t-F_t^\varepsilon)
  =-\int_0^\infty\varrho^\varepsilon(t+s)F_{-s}\dif s.
\]
The history norm gives
$\|F_{-s}\|_{H^m}\leq e^{s/2}s^{1/3}\|z_0\|_{\mathrm{hist}}$.
Together with~\eqref{eq:reg-mollifier-assumptions} and $\varepsilon\leq a$, this yields
\begin{equation}\label{eq:reg-remainder-mean}
  \|\E_{z_0}(F_t-F_t^\varepsilon)\|_{H^m}
  \leq C_0\varepsilon^{1/3}e^{-at/\varepsilon}\|z_0\|_{\mathrm{hist}}.
\end{equation}

Let $\widehat F$ have the same future as $F$ and an independent Brownian past, and set
$\widehat F^\varepsilon=\varrho^\varepsilon*\widehat F$.
The difference $\widehat F-\widehat F^\varepsilon$ is a continuous, stationary Gaussian process, and
$(\widehat F-\widehat F^\varepsilon)_{\varepsilon\,\cdot}$ has the same law as
$\sqrt{\varepsilon}(\widehat F-\widehat F^1)$.
Fernique's theorem in $C([0,1];H^m)$ therefore gives
\[
  \E\exp\left(\frac{c_0}{\varepsilon}
    \sup_{j\varepsilon\leq t\leq(j+1)\varepsilon}
      \|\widehat F_t-\widehat F_t^\varepsilon\|_{H^m}^2\right)\leq C
  \qquad(j\in\mathbb N_0).
\]
Since
\[
  F_t-F_t^\varepsilon-\E_{z_0}(F_t-F_t^\varepsilon)
  =\E\left[\widehat F_t-\widehat F_t^\varepsilon\mid(F_s)_{s\geq0}\right],
\]
conditional Jensen's inequality gives the same bound for
$F-F^\varepsilon-\E_{z_0}(F-F^\varepsilon)$ under $\mathbb P_{z_0}$.
Markov's inequality and a union bound, together with~\eqref{eq:reg-remainder-mean}, prove~\eqref{eq:reg-correction-tail}.
\end{proof}

The next estimate retains a strict contraction of the initial vorticity while allowing a fixed cost for the initial noise history. This cost will later be absorbed by choosing $\sigma$ small relative to $b$.

\begin{lemma}\label{lem:reg-vorticity-estimates}
For each $T\geq1$, there are $\rho,\eta_0,\varepsilon_0,C>0$ and $\vartheta\in(0,1)$ such that, for $0\leq\eta\leq\eta_0$ and $0\leq\varepsilon\leq\varepsilon_0$,
\begin{equation}\label{eq:Lyap-omega-part}
  \E_{z_0}\exp\big(\eta L_\rho(\omega_T^\varepsilon)\big)
  \leq C\exp\big(\eta\vartheta L_\rho(\omega_0)+C\eta\|z_0\|_{\mathrm{hist}}^2\big),
\end{equation}
and
\begin{equation}\label{eq:reg-energy-moment}
  \E_{z_0}\exp\big(\eta\mathcal E_T(\omega^\varepsilon)\big)
  \leq C\exp\big(C\eta(L_\rho(\omega_0)+\|z_0\|_{\mathrm{hist}}^2)\big),
\end{equation}
where
\begin{equation}\label{eq:reg-energy-functional}
  \mathcal E_T(\omega)\coloneq
  \sup_{0\leq t\leq T}\|\omega_t\|_{L^2}^2
  +\int_0^T\|\nabla\omega_t\|_{L^2}^2\dif t
  +\left(\sup_{0\leq t\leq T}\|\nabla\omega_t\|_{L^2}^2
       +\int_0^T\|\nabla^2\omega_t\|_{L^2}^2\dif t\right)^{1/3}.
\end{equation}
\end{lemma}
\begin{proof}
We treat $\varepsilon>0$ by splitting according to whether
$\|z_0\|_{\mathrm{hist}}^2\geq\varepsilon^{-1}$.
Throughout, $0<\varepsilon\leq\min\{1,a\}$.

\textbf{Step 1.}
We estimate exponential moments of
\[
  J_\varepsilon\coloneq
  \|z_0\|_{\mathrm{hist}}^2
  +\int_0^T
    \|\xi_t^\varepsilon-\E_{z_0}\xi_t^\varepsilon\|_{H^m}^2\dif t.
\]
The definition~\eqref{eq:F-eps-def}, the mollifier assumptions and the history norm give
\begin{equation}\label{eq:forcing-estimate-74}
\begin{aligned}
  \|\E_{z_0}\xi_t^\varepsilon\|_{H^m}
  &=\left\|\int_0^\infty(\varrho^\varepsilon)'(t+s)F_{-s}\dif s\right\|_{H^m}\\
  &\leq C\varepsilon^{-2/3}e^{-at/\varepsilon}\|z_0\|_{\mathrm{hist}},
\end{aligned}
\end{equation}
so
\begin{equation*}
  \left(\int_0^T \|\xi^\varepsilon_t\|_{H^m} \dif t\right)^2 \leq C J_\varepsilon.
\end{equation*}
Let $K^\varepsilon$ be the covariance operator of $\xi^\varepsilon$, viewed as a random variable on $L^2_t([0, T]; H^m)$.
Then $K^\varepsilon$ is trace class with
\begin{align*}
  \tr K^\varepsilon &= \E_{z_0} \int_0^T \|\xi^\varepsilon_t - \E_{z_0} \xi^\varepsilon_t\|_{H^m}^2 \dif t\\
                    &= \E_{z_0} \int_0^T \left\|\int_0^t \varrho^\varepsilon(t-s) \dif F_s\right\|_{H^m}^2 \dif t\\
                    &\leq C\E_{z_0} \int_0^T \int_0^t (\varrho^\varepsilon(t-s))^2 \dif s \dif t\\
                    &\leq C \varepsilon^{-1}.
\end{align*}
We estimate the largest eigenvalue of $K^\varepsilon$ by testing against $\varphi \in L^2_tH^m_x$:
\begin{align*}
  \E_{z_0} \langle \xi^\varepsilon - \E_{z_0} \xi^\varepsilon, \varphi \rangle^2
  &= \E_{z_0} \left(\int_0^T \int_0^t \varrho^\varepsilon(t-s) \langle \dif F_s, \varphi_t\rangle \, \dif t\right)^2\\
  &= \E_{z_0} \left(\sum_k \int_0^T \int_s^T \varrho^\varepsilon(t-s) \beta_k \langle e_k, \varphi_t \rangle \dif t \, \dif W^k_s\right)^2\\
  &\leq C \int_0^T \left(\int_s^T \varrho^\varepsilon(t-s) \|\varphi_t\|_{H^m} \dif t\right)^2 \, \dif s\\
  &\leq C \int_0^T \|\varphi_t\|_{H^m}^2 \dif t
\end{align*}
where we used It\^o's isometry, stochastic Fubini, and Young's convolution inequality.
So the largest eigenvalue of $K^\varepsilon$ is at most some constant $C$, which does not depend on $\varepsilon$, and therefore
\begin{equation}\label{eq:reg-crude-noise-moment}
  \E_{z_0}e^{\eta J_\varepsilon}
  \leq\exp\left(C\eta\left(\varepsilon^{-1}
                    +\|z_0\|_{\mathrm{hist}}^2\right)\right),
  \qquad 0\leq\eta\leq (4C)^{-1}.
\end{equation}

\textbf{Step 2.} We prove the lemma when $\|z_0\|_{\mathrm{hist}}$ is large.

Test~\eqref{eq:2dnse-smooth-noise} against $e^{\nu t}\omega^\varepsilon$ and $-e^{\nu t / 2}\Delta\omega^\varepsilon$.
For every $\zeta>0$ and $0\leq t\leq T$, this gives
\begin{align}
  &e^{\nu t}\|\omega_t^\varepsilon\|_{L^2}^2
       +\nu\int_0^t e^{\nu s}\|\nabla\omega_s^\varepsilon\|_{L^2}^2\dif s
  \leq(1+\zeta)\|\omega_0\|_{L^2}^2
       +C_{\zeta,T}J_\varepsilon,\label{eq:reg-deterministic-energy}\\
  &e^{\nu t/2}\|\nabla\omega_t^\varepsilon\|_{L^2}^2
       +\frac{\nu}{2}\int_0^t e^{\nu s/2}\|\nabla^2\omega_s^\varepsilon\|_{L^2}^2\dif s
       \label{eq:reg-deterministic-energy-2}\\
  &\qquad\leq(1+\zeta)\|\nabla\omega_0\|_{L^2}^2
       +C_{\zeta,T}\int_0^t e^{\nu s/2}\|\omega_s^\varepsilon\|_{L^2}^6\dif s
       +C_{\zeta,T}J_\varepsilon.\nonumber
\end{align}
Choosing $\zeta>0$ and then $\rho>0$ sufficiently small gives some $\vartheta<1$ such that
\begin{equation}\label{eq:reg-crude-pathwise}
  L_\rho(\omega_T^\varepsilon)
  \leq\vartheta L_\rho(\omega_0)+C(1+J_\varepsilon).
\end{equation}
Meanwhile, any choice of $\zeta, \rho > 0$ yields
\begin{equation}\label{eq:reg-crude-pathwise-2}
  \mathcal E_T(\omega^\varepsilon)
  \leq C\big(L_\rho(\omega_0)+1+J_\varepsilon\big).
\end{equation}

Suppose $\|z_0\|_{\mathrm{hist}}^2\geq\varepsilon^{-1}$.
Then~\eqref{eq:reg-crude-noise-moment} gives
\[
  \E_{z_0}e^{\eta J_\varepsilon}
  \leq e^{C_T\eta\|z_0\|_{\mathrm{hist}}^2}.
\]
Combining this with~\eqref{eq:reg-crude-pathwise}--\eqref{eq:reg-crude-pathwise-2} proves both assertions for sufficiently small $\eta$.

\textbf{Step 3.} We prove the lemma when $\|z_0\|_{\mathrm{hist}}$ is not too large.

Suppose $\|z_0\|_{\mathrm{hist}}^2<\varepsilon^{-1}$.
Define the remainder $R^\varepsilon\coloneq F-F^\varepsilon$ and its deterministic history bound
\[
  q_\varepsilon(t)\coloneq
  C_0\varepsilon^{1/3}e^{-at/\varepsilon}\|z_0\|_{\mathrm{hist}}.
\]
The assumption on the history gives
\begin{equation}\label{eq:reg-small-history-bounds}
  \int_0^Tq_\varepsilon(t)\dif t\leq C\varepsilon^{5/6},
  \qquad
  \sup_{[0, T]}q_\varepsilon^2
       +\int_0^Tq_\varepsilon(t)^4\dif t
  \leq C\|z_0\|_{\mathrm{hist}}^2.
\end{equation}
Fix $\theta>0$ sufficiently small depending on $\nu$, and define
\[
  \mathcal A_\varepsilon(z_0)=
  \left\{\|R_t^\varepsilon\|_{H^m}\leq\theta+q_\varepsilon(t)
         \text{ for every }0\leq t\leq T\right\}.
\]
By Lemma~\ref{lem:reg-noise-estimates}, for some $c_\theta>0$,
\begin{equation}\label{eq:reg-exceptional-probability}
  \mathbb P_{z_0}(\mathcal A_\varepsilon(z_0)^c)
  \leq C e^{-c_\theta/\varepsilon}.
\end{equation}

Let $v^\varepsilon=\omega^\varepsilon+R^\varepsilon$.
Then $v^\varepsilon$ almost solves the same equation as $\omega^0$:
\begin{equation}\label{eq:reg-compensated-vorticity}
\begin{aligned}
  \dif v^\varepsilon
  &=\big[\nu\Delta(v^\varepsilon-R^\varepsilon)\\
  &\qquad-\big(\nabla^\perp\Delta^{-1}(v^\varepsilon-R^\varepsilon)\big)
       \cdot\nabla(v^\varepsilon-R^\varepsilon)\big]\dif t+\dif F_t,\\
  v_0^\varepsilon&=\omega_0+R_0^\varepsilon.
\end{aligned}
\end{equation}

\textbf{Step 3a.} We prove~\eqref{eq:Lyap-omega-part} and~\eqref{eq:reg-energy-moment} with $\omega^\varepsilon$ replaced by the corrected vorticity $v^\varepsilon$.

Define the stopping time $t_1 \coloneq T \land \inf \left\{t \in [0, T] \mid \|R_t^\varepsilon\|_{H^m}>\theta+q_\varepsilon(t)\right\}$.
Before time $t_1$, testing against $v^\varepsilon_t$ with It\^o's formula gives a constant $C_*>0$ such that
\begin{align}\label{eq:reg-stopped-energy}
  \dif\|v_t^\varepsilon\|_{L^2}^2
       +\nu\|\nabla v_t^\varepsilon\|_{L^2}^2\dif t
  &\leq C_*q_\varepsilon(t)\|v_t^\varepsilon\|_{L^2}^2\dif t\nonumber\\
  &\quad+C\big(1+q_\varepsilon(t)^4\big)\dif t
       +2\langle v_t^\varepsilon,\dif F_t\rangle,\\
  \dif\|\nabla v_t^\varepsilon\|_{L^2}^2
       +\nu\|\nabla^2v_t^\varepsilon\|_{L^2}^2\dif t
  &\leq C_*q_\varepsilon(t)\|\nabla v_t^\varepsilon\|_{L^2}^2\dif t\nonumber\\
  &\quad+C\big(1+\|v_t^\varepsilon\|_{L^2}^6
                 +q_\varepsilon(t)^4\big)\dif t
       +2\langle\nabla v_t^\varepsilon,\nabla\dif F_t\rangle.\nonumber
\end{align}

For $0<\lambda\leq\nu/2$, define the integrating factor
\[
  f_\varepsilon(t)\coloneq
  \exp\left(\lambda t-C_*\int_0^tq_\varepsilon(s)\dif s\right).
\]
Choose $\varepsilon_0>0$ sufficiently small that
$C_*\int_0^Tq_\varepsilon(s)\dif s\leq\lambda/2$
whenever $\varepsilon\leq\varepsilon_0$ and
$\|z_0\|_{\mathrm{hist}}^2<\varepsilon^{-1}$.
Then
\begin{equation}\label{eq:reg-integrating-factor}
  f_\varepsilon'=(\lambda-C_*q_\varepsilon)f_\varepsilon,\qquad
  e^{-\lambda/2}\leq f_\varepsilon(t)\leq e^{\lambda T},\qquad
  f_\varepsilon(T)\geq e^{\lambda T/2}>1.
\end{equation}
Multiplying~\eqref{eq:reg-stopped-energy} by $f_\varepsilon$ and using Poincar\'e's inequality gives, for $0\leq t\leq t_1$,
\begin{align}
  f_\varepsilon(t)\|v_t^\varepsilon\|_{L^2}^2 +\frac{\nu}{2}\int_0^t f_\varepsilon(s)\|\nabla v_s^\varepsilon\|_{L^2}^2\dif s
  &\leq \|v^\varepsilon_0\|_{L^2}^2 + C\int_0^t f_\varepsilon(s)(1 + q_\varepsilon(s)^4) \dif s + M_0(t),\label{eq:final-martingale-1}\\
  f_\varepsilon(t)\|\nabla v_t^\varepsilon\|_{L^2}^2 +\frac{\nu}{2}\int_0^t f_\varepsilon(s)\|\nabla^2 v_s^\varepsilon\|_{L^2}^2\dif s
  &\leq \|\nabla v^\varepsilon_0\|_{L^2}^2 + C\int_0^t f_\varepsilon(s)(1 + \|v^\varepsilon_s\|_{L^2}^6 + q_\varepsilon(s)^4) \dif s + M_1(t)\label{eq:final-martingale-2},
\end{align}
where, for $0\leq t\leq T$,
\begin{align*}
  M_0(t)&\coloneq
    2\int_0^{t\land t_1}
      f_\varepsilon(s)\langle v_s^\varepsilon,\dif F_s\rangle_{L^2},\\
  M_1(t)&\coloneq
    2\int_0^{t\land t_1}
      f_\varepsilon(s)\langle\nabla v_s^\varepsilon,\nabla\dif F_s\rangle_{L^2}.
\end{align*}

These are square-integrable martingales by Steps~1--2 and the bound on $R^\varepsilon$ up to $t_1$.
For every $c>0$, the processes $\exp(cM_i(t)-[cM_i]_t)$ are nonnegative supermartingales starting at one, and the quadratic variations are bounded by
\begin{align*}
  [M_0]_t
  &\leq C\int_0^{t\land t_1}
      f_\varepsilon(s)^2\|v_s^\varepsilon\|_{L^2}^2\dif s
  \leq C\int_0^{t\land t_1}
      f_\varepsilon(s)\|v_s^\varepsilon\|_{L^2}^2\dif s,\\
  [M_1]_t
  &\leq C\int_0^{t\land t_1}
      f_\varepsilon(s)^2\|\nabla v_s^\varepsilon\|_{L^2}^2\dif s
  \leq C\int_0^{t\land t_1}
      f_\varepsilon(s)\|\nabla v_s^\varepsilon\|_{L^2}^2\dif s.
\end{align*}
Choose $c>0$ sufficiently small that Poincar\'e's inequality absorbs $c[M_i]_t$ into half of the corresponding dissipative term in~\eqref{eq:final-martingale-1}--\eqref{eq:final-martingale-2}.
Doob's maximal inequality then gives exponential tails for $\sup_{t\leq T}(M_i(t)-c[M_i]_t)$.
Take a cube root in~\eqref{eq:final-martingale-2} and use~\eqref{eq:final-martingale-1} to control the integral of $\|v_s^\varepsilon\|_{L^2}^6$.
Since $t_1=T$ on $\mathcal A_\varepsilon(z_0)$, the bounds~\eqref{eq:reg-small-history-bounds} and~\eqref{eq:reg-integrating-factor} give, after decreasing $\rho,\eta_0>0$ if necessary,
\begin{align*}
  \E_{z_0}\left[\mathbf1_{\mathcal A_\varepsilon(z_0)}
    \exp\left({\eta e^{\lambda T/6}L_\rho(v_T^\varepsilon)}\right)\right]
  &\leq C\exp\big(\eta L_\rho(v_0^\varepsilon)+C\eta\|z_0\|_{\mathrm{hist}}^2\big),\\
  \E_{z_0}\left[\mathbf1_{\mathcal A_\varepsilon(z_0)}
    \exp\left({\eta\mathcal E_T(v^\varepsilon)}\right)\right]
  &\leq C\exp\big(C\eta(L_\rho(v_0^\varepsilon)+\|z_0\|_{\mathrm{hist}}^2)\big)
\end{align*}
for $0\leq\eta\leq\eta_0$.

\textbf{Step 3b.} We transfer the previous estimates to $\omega^\varepsilon$.

Work in the event $\mathcal{A}_\varepsilon(z_0)$. For every $\zeta>0$, the bound $\|R_t^\varepsilon\|_{H^m}\leq\theta+q_\varepsilon(t)$ on $\mathcal A_\varepsilon(z_0)$ gives
\begin{align*}
  L_\rho(\omega_T^\varepsilon)
  &\leq(1+\zeta)L_\rho(v_T^\varepsilon)+C_\zeta(1+\|z_0\|_{\mathrm{hist}}^2),\\
  L_\rho(v_0^\varepsilon)
  &\leq(1+\zeta)L_\rho(\omega_0)+C_\zeta(1+\|z_0\|_{\mathrm{hist}}^2),\\
  \mathcal E_T(\omega^\varepsilon)
  &\leq C\mathcal E_T(v^\varepsilon)+C(1+\|z_0\|_{\mathrm{hist}}^2).
\end{align*}
We conclude by choosing $\zeta > 0$ sufficiently small.

In the complement of $\mathcal{A}_\varepsilon(z_0)$, Cauchy--Schwarz and~\eqref{eq:reg-crude-pathwise},~\eqref{eq:reg-crude-noise-moment}, and~\eqref{eq:reg-exceptional-probability} give
\[
  \E_{z_0}\big[\mathbf1_{\mathcal A_\varepsilon(z_0)^c}
    e^{\eta L_\rho(\omega_T^\varepsilon)}\big]
  \leq C e^{\eta\vartheta L_\rho(\omega_0)+C\eta\|z_0\|_{\mathrm{hist}}^2}
    \exp\big((C\eta-c_\theta/2)/\varepsilon\big),
\]
which (up to decreasing $\eta_0$ further) proves~\eqref{eq:Lyap-omega-part}.
Conclusion~\eqref{eq:reg-energy-moment} follows similarly, using~\eqref{eq:reg-crude-pathwise-2}.

The case $\varepsilon = 0$ is strictly easier, since $R^0 = 0$.
\end{proof}

\begin{proof}[Proof of Proposition~\ref{prop:reg-V-NSE}]
Splitting the supremum in the history norm according to whether $T+s\leq0$ or $T+s\geq0$ gives
\[
  \|z_T^\varepsilon\|_{\mathrm{hist}}^2
  \leq2e^{-T}\|z_0\|_{\mathrm{hist}}^2
       +2\sup_{0\leq s<t\leq T}
         \frac{\|F_t-F_s\|_{H^m}^2}{|t-s|^{2/3}}.
\]
The supremum is the square of the Brownian $1/3$-H\"older seminorm on $[0,T]$, which has an exponential moment at a sufficiently small exponent. Therefore, for small $a_1>0$ we have
\begin{equation}\label{eq:Lyap-F-part}
  \E_{z_0}e^{\eta \|z_T^\varepsilon\|_{\mathrm{hist}}^2}\leq C e^{2\eta e^{-T}\|z_0\|_{\mathrm{hist}}^2},
  \qquad 0\leq\eta\leq a_1.
\end{equation}

Take $\rho,\vartheta,\eta_0$ from Lemma~\ref{lem:reg-vorticity-estimates} and choose
\[
  1<\alpha<\min\{\vartheta^{-1},e/2\}.
\]
First choose $b>0$ so that $2\alpha b\leq a_1$, and then choose $\sigma>0$ so that
\[
  2\alpha\sigma\leq\eta_0,
  \qquad \alpha C\sigma+2\alpha b e^{-T}\leq b,
\]
where $C$ is the history coefficient in~\eqref{eq:Lyap-omega-part}. By Cauchy--Schwarz,
\begin{align*}
  \E_{z_0}e^{\alpha V_{\sigma,\rho,b}(z_T^\varepsilon)}
  &\leq\big(\E_{z_0}e^{2\alpha\sigma L_\rho(\omega_T^\varepsilon)}\big)^{1/2}
        \big(\E_{z_0}e^{2\alpha b\|z_T^\varepsilon\|_{\mathrm{hist}}^2}\big)^{1/2}\\
  &\leq C\exp\big(\alpha\sigma\vartheta L_\rho(\omega_0)
                    +(\alpha C\sigma+2\alpha b e^{-T})\|z_0\|_{\mathrm{hist}}^2\big)\\
  &\leq C\exp\big(V_{\sigma,\rho,b}(z_0)\big),
\end{align*}
proving~\eqref{eq:reg-V-contraction-NSE}.

To prove the remaining estimates, choose $0<r<s<1$ so small that
\[
  \alpha(1/2+s)<1.
\]
Sobolev embedding, interpolation, and H\"older's inequality in time give
\begin{align*}
  \int_0^T\|\nabla u_t^\varepsilon\|_{C^r}^{\alpha}\dif t
  &\leq C\int_0^T\|\omega_t^\varepsilon\|_{H^1}^{\alpha(1-s)}
                         \|\omega_t^\varepsilon\|_{H^2}^{\alpha s}\dif t\\
  &\leq C\left(\int_0^T\|\omega_t^\varepsilon\|_{H^1}^2\dif t\right)^{\alpha(1-s)/2}
          \left(\int_0^T\|\omega_t^\varepsilon\|_{H^2}^2\dif t\right)^{\alpha s/2}\\
  &\leq C\mathcal E_T(\omega^\varepsilon)^{\alpha(1/2+s)}.
\end{align*}
For every $A>0$ and every sufficiently small $\eta>0$, Young's inequality therefore gives
\[
  A\int_0^T\|\nabla u_t^\varepsilon\|_{C^r}^{\alpha}\dif t
  \leq\eta\mathcal E_T(\omega^\varepsilon)+C_{A,\eta}.
\]
Using~\eqref{eq:reg-energy-moment} and choosing $\eta$ so that $C\eta\leq\delta\min\{\sigma,b\}$ proves~\eqref{eq:reg-V-bounded-NSE}, and in fact the same bound with $A$ in place of $r$. Finally, the flow estimate in Lemma~\ref{lem:velocityToFlowMap} bounds $\sup_{t\leq T}\BiLip(\varphi_t^\varepsilon)^p$ by a constant multiple of $\exp(C_p\int_0^T\|\nabla u_t^\varepsilon\|_{L^\infty}\dif t)$. Since $\alpha>1$, the preceding bound with an appropriate $A$ also proves~\eqref{eq:reg-phi-bounded-NSE}.
\end{proof}

\subsection{Stability and mixing}

We first record the white-noise mixing estimate in a form that allows $H^1$ initial vorticity.

\begin{lemma}\label{lem:reg-white-mixing}
There are $\delta\in(0,1)$ and $\mu>0$ such that, for every $p\in(1,2)$, $\rho>0$ and $\zeta>0$, there is $C>0$ with
\begin{equation}\label{eq:reg-white-mixing}
  \left(\E_{\omega_0}\|\theta_n^0\|_{H^{-\delta}}^p\right)^{1/p}
  \leq C e^{-\mu n}e^{\zeta L_\rho(\omega_0)}\|\theta_0\|_{H^{-\delta}}
  \qquad(n\in\mathbb N)
\end{equation}
for the white-noise equation, every $\omega_0\in H^1_0$ and every $\theta_0\in H^{-\delta}_0$.
\end{lemma}
\begin{proof}
The velocity noise has Fourier coefficients comparable to $|k|^{-(\ell+1)}$, and
\[
  \ell+1>5,\qquad \ell-1<m+1<\ell,\qquad m+1>4.
\]
Thus~\cite[Theorem~1.8 and Section~2.2]{BFPS24} applies with velocity phase space $H^{m+1}$, equivalently vorticity phase space $H^m$. It gives $\delta,\mu>0$ and $M>0$ such that, for sufficiently small $\eta>0$,
\begin{equation}\label{eq:reg-white-smooth-mixing}
  \left(\E_{\omega_0}\|\theta_n^0\|_{H^{-\delta}}^p\right)^{1/p}
  \leq C_\eta e^{-\mu n}(1+\|\omega_0\|_{H^m})^M
           e^{\eta\|\omega_0\|_{L^2}^2}\|\theta_0\|_{H^{-\delta}}.
\end{equation}

First, observe that~\eqref{eq:reg-white-mixing} is satisfied when $n=1$ by choosing $C > 0$ large enough and applying~\eqref{eq:reg-phi-bounded-NSE}.

Then test~\eqref{eq:2dnse-smooth-noise} against $t^m(-\Delta)^m\omega^0$ to obtain
\begin{equation}\label{eq:reg-white-smoothing}
  \E_{\omega_0}(1+\|\omega_1^0\|_{H^m})^N
  \leq C_\zeta e^{\zeta L_\rho(\omega_0)}
\end{equation}
for arbitrarily small $\zeta > 0$, and conclude~\eqref{eq:reg-white-mixing} for all $n \geq 2$ by composing with~\eqref{eq:reg-white-smooth-mixing} from time $1$ to $n$.
\end{proof}

We now fix $V=V_{\sigma,\rho,b}$, $\alpha>1$, $K>1$ and $\varepsilon_0>0$ given by Proposition~\ref{prop:reg-V-NSE} at $T=1$. Thus
\begin{equation}\label{eq:reg-fixed-drift}
  \E_{z_0} e^{\alpha V(z_1^\varepsilon)}\leq K e^{V(z_0)}
  \qquad(0\leq\varepsilon\leq\varepsilon_0).
\end{equation}
By the Markov property and Jensen's inequality, for $j\in \Z_{\geq 0}$ and $0<c\leq\alpha$, we inductively have
\begin{equation}\label{eq:reg-iterated-drift}
  \E_{z_0} e^{cV(z_j^\varepsilon)}
  \leq K^{c/(\alpha-1)}e^{c\alpha^{-j}V({z_0})}.
\end{equation}
Also, for each fixed $S\geq1$ and $p,\zeta>0$ and all sufficiently small $\varepsilon\geq0$, Proposition~\ref{prop:reg-V-NSE} applied at each integer time in $[0, \ceil{S}-1]$ with~\eqref{eq:reg-iterated-drift} gives
\begin{equation}\label{eq:reg-fixed-weight-flow}
  \E_{z_0}\sup_{0\leq t\leq S}\BiLip(\varphi_t^\varepsilon)^p
  \leq C_{S,p,\zeta}e^{\zeta V({z_0})}.
\end{equation}

Next, we compare the white-in-time and regularized systems. Given $z_0\in Z$ and $\varepsilon>0$, let $\omega^\varepsilon$ and $\omega^0$ start from the same $\omega_0$ and use the same future Brownian increments of $F$. Their flow maps from time zero are $\varphi_t^\varepsilon$ and $\varphi_t^0$.

\begin{lemma}\label{lem:reg-comparison}
For every fixed $S\geq1$ and every $p>0$ and $M\geq0$,
\begin{equation}\label{eq:reg-weighted-flow-convergence}
  \lim_{\varepsilon \to 0}\sup_{z_0\in Z}e^{-pV(z_0)}
  \E_{z_0}\left[(\BiLip(\varphi_S^\varepsilon)+\BiLip(\varphi_S^0))^M
      d^1(\varphi_S^\varepsilon,\varphi_S^0)\right]=0.
\end{equation}
\end{lemma}
\begin{proof}
Fix $S \geq 1$, $p > 0$, $M \geq 0$, and $z_0 \in Z$. By Lemma~\ref{lem:velocityToFlowMap} and interpolation with~\eqref{eq:reg-energy-moment} it suffices to prove that
\begin{equation*}
  \E_{z_0}\left[\int_0^S \|\omega^\varepsilon_t - \omega^0_t\|_{L^2}^c \dif t\right]
  = e^{pV(z_0)}o_{\varepsilon}(1),
\end{equation*}
for some small constant $c > 0$, where the $o_{\varepsilon}(1)$ term does not depend on $z_0$.
We first prove convergence in probability uniformly on each set $\{V(z_0)\leq R\}$. Write $R^\varepsilon=F-F^\varepsilon$ and put
\[
  e_t^\varepsilon=\omega_t^\varepsilon+R_t^\varepsilon-\omega_t^0,
  \qquad 0\leq t\leq S.
\]
Then $e_0^\varepsilon=R_0^\varepsilon$, and the Brownian terms cancel, giving
\begin{equation*}
  \partial_t e^\varepsilon
  +u^\varepsilon\cdot\nabla(e^\varepsilon-R^\varepsilon)
  +\big(\nabla^\perp\Delta^{-1}(e^\varepsilon-R^\varepsilon)\big)
  \cdot\nabla\omega^0
  -\nu\Delta(e^\varepsilon-R^\varepsilon)
  =0.
\end{equation*}
Testing this equation against $e^\varepsilon$ gives
\begin{equation*}
  \frac{\dif}{\dif t}\|e_t^\varepsilon\|_{L^2}^2+\nu\|\nabla e_t^\varepsilon\|_{L^2}^2
  \leq C(1+\|\omega_t^0\|_{L^4}^2)\|e_t^\varepsilon\|_{L^2}^2
  +C\|R_t^\varepsilon\|_{H^m}^2
       (1+\|\omega_t^\varepsilon\|_{L^2}^2+\|\omega_t^0\|_{L^2}^2).
\end{equation*}
Then Gr\"onwall's inequality with initial condition $\|e^\varepsilon_0\|_{L^2}^2 = \|R^\varepsilon_0\|_{L^2}^2$ gives
\begin{equation*}
  \|e^\varepsilon_t\|_{L^2}^2 + \nu \int_0^t \|\nabla e^\varepsilon_s\|_{L^2}^2 \dif s \leq C\left(\sup_{0 \leq s \leq t}\|R^\varepsilon_s\|_{H^m}^2\right)(1 + \mathcal{E}_t(\omega^0) + \mathcal{E}_t(\omega^\varepsilon))\exp(C\mathcal{E}_t(\omega^0)),
\end{equation*}
and therefore
\begin{equation}\label{eq:almost-done}
  \|\omega^\varepsilon_t - \omega^0_t\|_{L^2}^2 \leq C\left(\sup_{0 \leq s \leq t}\|R^\varepsilon_s\|_{H^m}^2\right)(1 + \mathcal{E}_t(\omega^0) + \mathcal{E}_t(\omega^\varepsilon))\exp(C\mathcal{E}_t(\omega^0)).
\end{equation}
Choose $c$ small enough so that $Cc \ll \eta_0 \land p\sigma \land pb$ from Lemma~\ref{lem:reg-vorticity-estimates},
and raise~\eqref{eq:almost-done} to the power $c$, integrate from $0$ to $S$ in time, and take an expectation to conclude.
\end{proof}

\begin{proof}[Proof of Theorem~\ref{thm:reg-noise-mixing}]
  Assume without loss of generality that $\omega_0 \in H^2_0(\T^2)$ for compatibility with the state space $Z$. Indeed, the conclusion of the theorem depends only on the $H^1$ norm and therefore we recover the same result for $\omega_0 \in H^1_0(\T^2)$ by approximation.

Fix $p\in(1,2)$ and $\delta,\mu$ as in Lemma~\ref{lem:reg-white-mixing}, and retain the fixed $V,\alpha,K$ from~\eqref{eq:reg-fixed-drift}. Put
\[
  \overline K=K^{\alpha/(\alpha-1)}.
\]
For every integer $N\geq1$,~\eqref{eq:reg-iterated-drift} implies
\begin{equation}\label{eq:reg-block-drift}
  \E_{z_0} e^{\alpha V(z_N^\varepsilon)}\leq\overline K e^{V(z_0)}.
\end{equation}
In particular, this drift constant is independent of $N$.

Choose an integer $N\geq2$ sufficiently large that the contraction factor
\[
  \kappa_N\coloneq C_0e^{-\mu N}
  <\frac14\overline K^{-1/(\alpha-1)}
\]
is compatible with the Lyapunov function as in Assumption~\ref{ass:super_lyapunov}.

Consider two transition kernels on $Z$ over time $N$: the reference kernel $P^0$ uses the white-noise equation on $[0,N]$, while $P^\varepsilon$ uses the regularized equation throughout.
Let $\varphi^\varepsilon_z \in \Diff^1(\T^2)$ be the flow map from time $-N$ to $0$ of the velocity $u^\varepsilon_t \coloneq \nabla^\perp \Delta^{-1} \omega^\varepsilon_t$ (where we allow $\varepsilon = 0$).
Observe that indeed $\varphi^\varepsilon_z \in \Diff^1(\T^2)$ almost surely by Lemma~\ref{lem:reg-comparison}.

By Lemma~\ref{lem:reg-white-mixing}, the kernel $P^0$ satisfies Assumption~\ref{ass:negative_regularity_mixing} with $V,p,\delta$ and $\kappa_N$. Lemma~\ref{lem:reg-comparison}, with $S=N$ and $M=p(1+\delta)+1$, verifies~\eqref{eq:smallness_condition_Wasmodified} for small $\varepsilon > 0$. Theorem~\ref{thm:mainStability} therefore gives a contraction factor $\kappa > 0$ for $P^\varepsilon$ with
\[
  \kappa\overline K^{1/(\alpha-1)}\leq\frac12<1
\]
for every $0<\varepsilon\leq\varepsilon_0$, after fixing $\varepsilon_0>0$ sufficiently small. Together with~\eqref{eq:reg-block-drift}, this verifies Assumption~\ref{ass:super_lyapunov} for the regularized system.

Theorem~\ref{thm:shortTimeToLong} now applies to the time-$N$ chain. In particular, choose
\[
  0<q_0<\frac{p(\alpha-1)}{p+\alpha-1},
  \qquad \chi=\frac{\alpha}{\alpha-1}.
\]
There are $C,\gamma_0>0$, independent of $\varepsilon\leq\varepsilon_0$, such that
\begin{equation}\label{eq:reg-conditional-block-mixing}
  \left(\E_{z_0}\|\theta_{jN}^\varepsilon\|_{H^{-\delta}}^{q_0}\right)^{1/q_0}
  \leq C e^{-\gamma_0j}e^{\chi V(z_0)}\|\theta_0\|_{H^{-\delta}}
  \qquad(j\in\mathbb N_0).
\end{equation}

To include intermediate times, take $0<q<q_0/2$ and write $t=jN+s$, with $0\leq s<N$. By Lemma~\ref{lem:HDeltaPerturb} and Cauchy--Schwarz,
\begin{align*}
  \left(\E_{z_0}\|\theta_t^\varepsilon\|_{H^{-\delta}}^q\right)^{1/q}
  &\leq\left(\E_{z_0}\|\theta_{jN}^\varepsilon\|_{H^{-\delta}}^{2q}\right)^{1/(2q)}
       \left(\E_{z_0}\BiLip(\varphi_{jN,t}^\varepsilon)^{2q(1+\delta)}\right)^{1/(2q)}.
\end{align*}
Conditioning the second factor at time $jN$, and using~\eqref{eq:reg-fixed-weight-flow} with initial weight $qV$, followed by~\eqref{eq:reg-iterated-drift}, bounds it by $Ce^{V(z_0)/2}$. The first factor is bounded by~\eqref{eq:reg-conditional-block-mixing}, since $2q<q_0$. Thus, after changing $C,\gamma>0$,
\begin{equation}\label{eq:reg-conditional-mixing}
  \left(\E_{z_0}\|\theta_t^\varepsilon\|_{H^{-\delta}}^q\right)^{1/q}
  \leq C e^{-\gamma t}e^{(\chi+1/2)V(z_0)}\|\theta_0\|_{H^{-\delta}}
  \qquad(t\geq0).
\end{equation}

Finally, averaging over the Brownian initial condition for the noise, we have $\E e^{b_0\|z_0\|_{\mathrm{hist}}^2}<\infty$ for some $b_0>0$. Decrease $q>0$ once more so that $q(\chi+1/2)b<b_0$ and conclude by taking the expectation of~\eqref{eq:reg-conditional-mixing} over initial noise histories.
\end{proof}

\appendix

\section{A corollary of Doob's maximal inequality}
\begin{lemma}\label{lemma:doobvariant}
Let $(\Omega,\mathscr F,\{\mathscr F_t\}_{0\leq t <\infty},\mathbb P)$ be a filtered probability space and $\{W_t\}_{0\leq t < \infty}$ be a $\mathscr F_t$-adapted Brownian motion.
Let $\{M_t\}_{0\leq t \leq T}$ be the martingale given by the stochastic integral
\[
M_t = M_0 + \int_0^t X_s \mathrm dW_s,
\]
where $M_0$ is $\mathscr F_0$-measurable and $X_s$ is progressively measurable. Let $[M]_t = \int_0^t X_s^2 ds$ be its quadratic variation on $[0,t]$. Let $ p > 1$ and assume the Novikov condition
\[
\mathbf E\left[\exp\left(\frac{p^2}{2} [M]_T\right)\right] <\infty.
\]
Then, we have
\[
\mathbf E\left[\sup_{t \in [0,T]}\exp\left(M_t - \frac{p}{2} [M]_t\right)\right] \leq \frac{p}{p-1}\mathbf E\left[\exp\left(M_0\right)\right].
\]
\end{lemma}
\begin{proof}
  Without loss of generality, we assume that $M_0$ is deterministic.
  Under the stated Novikov condition, $\exp\left(p(M_t- M_0) - \frac{p^2}{2} [M]_t\right)$ is a martingale on $[0,T]$ with mean 1~\cite[Chapter 3, Section 5]{karatzas1991brownian}. Hence, by Doob's maximal inequality, for all $\alpha >0$ we have
  \[
    \mathbb P\left[\sup_{t \in [0,T]}\exp\left(p(M_t-M_0) - \frac{p^2}{2} [M]_t\right) > \alpha\right]\leq\min\left\{1,\frac{1}{\alpha}\right\}.
  \]
  Thus,
  \begin{align*}
    \mathbf E\left[\sup_{t \in [0,T]}\exp\left((M_t- M_0) - \frac{p}{2} [M]_t\right)\right] &= \int_0^\infty \mathbb P\left[\sup_{t \in [0,T]}\exp\left(p(M_t-M_0) - \frac{p^2}{2} [M]_t\right)> \alpha^p\right] d\alpha \\
                                                                                            &\leq \int_0^\infty \min\left\{1,\frac{1}{\alpha^p}\right\}d\alpha = \frac{p}{p-1}.
  \end{align*}
\end{proof}

\bibliographystyle{alpha}
\bibliography{references.bib}
%

\end{document}